\documentclass[12pt, colorinlistoftodos]{amsart}
\usepackage{a4wide}
\usepackage{amssymb}
\usepackage{amsthm}
\usepackage{amsmath}
\usepackage{amscd}
\usepackage{verbatim}
\usepackage{bm}
\usepackage[mathscr]{eucal}
\usepackage[all]{xy}
\usepackage{todonotes}
\usepackage[pagebackref=true]{hyperref}
\usepackage{dsfont, stmaryrd}
\usepackage{mathrsfs}
\usepackage{mathtools}
\numberwithin{equation}{section}

\theoremstyle{plain}
\newtheorem{theorem}{Theorem}[section]

\newtheorem{lemma}[theorem]{Lemma}
\newtheorem{proposition}[theorem]{Proposition}

\theoremstyle{definition}
\newtheorem{definition}[theorem]{Definition}
\newtheorem{remark}[theorem]{Remark}

\theoremstyle{remark}

\newcommand{\pmat}[4]{\begin{pmatrix}
                 #1 & #2\\
                 #3 & #4
\end{pmatrix}}
\newcommand{\smat}[4]{\left(\begin{smallmatrix}
                 #1 & #2\\
                 #3 & #4
\end{smallmatrix}\right)}

\newcommand{\lp}{\left (}
\newcommand{\rp}{\right )}

\newcommand{\Fc}{{\mathcal{F}}}
\newcommand{\Hc}{{\mathcal{H}}}
\newcommand{\Nc}{{\mathcal{N}}}
\newcommand{\Oc}{{\mathcal{O}}}
\newcommand{\Sc}{{\mathcal{S}}}

\newcommand{\Zb}{\mathbb{Z}}
\newcommand{\Fb}{\mathbb{F}}
\newcommand{\Qb}{\mathbb{Q}}
\newcommand{\SO}{{\mathrm{SO}}}
\newcommand{\GL}{{\mathrm{GL}}}
\newcommand{\sgn}{{\mathrm{sgn}}}

\newcommand{\nf}{\mathfrak{n}}

\newcommand{\df}{\mathfrak{d}}
\newcommand{\ef}{\mathfrak{e}}

\newcommand{\wf}{\mathfrak{w}}

\newcommand{\uf}{\mathfrak{u}}
\newcommand{\vf}{\mathfrak{v}}

\newcommand{\Nf}{\mathfrak{N}}

\newcommand{\Nm}{{\mathrm{Nm}}}
\newcommand{\Ab}{\mathbb{A}}
\newcommand{\Pb}{\mathbb{P}}

\newcommand{\Nb}{\mathbb{N}}
\newcommand{\Hb}{\mathbb{H}}
\newcommand{\Rb}{\mathbb{R}}
\newcommand{\Cb}{\mathbb{C}}

\newcommand{\ebf}{{\mathbf{e}}}

\newcommand{\tr}{\operatorname{Tr}}

\newcommand{\GSpin}{\operatorname{GSpin}}

\newcommand{\Gal}{\operatorname{Gal}}

\newcommand{\kro}[2]{\left( \tfrac{#1}{#2} \right)}
\newcommand{\Lo}{L}

\newcommand{\vol}{\mathrm{vol}}

\newcommand{\Fr}{\mathcal{F}}
\newcommand{\Cc}{\mathcal{C}}
\newcommand{\SL}{\mathrm{SL}}
\renewcommand{\k}{{E}}
\newcommand{\K}{\mathrm{K}^1}
\newcommand{\cha}{\mathrm{Char}}
\renewcommand{\aa}{d}

\newcommand{\tGoo}{{\tilde{G}_0}}
\newcommand{\tGo}{{\tilde{G}}}
\newcommand{\Go}{G}
\newcommand{\Goo}{G_0}

\newcommand{\dfn}{d}
\newcommand{\Wb}{\mathbb{W}}
\newcommand{\Wc}{\mathcal{W}}
\newcommand{\alb}{{\boldsymbol{ \alpha}}}
\renewcommand{\Mc}{\mathcal{M}}

\newcommand{\Ec}{\mathcal{E}}
\newcommand{\tv}{{\tilde\varphi}}

\newcommand{\tM}{{\tilde{M}}}
\newcommand{\tQ}{{\tilde{Q}}}
\newcommand{\Vo}{{V^\circ}}
\newcommand{\Ho}{H^\circ}
\newcommand{\vv}{\mathrm{vv}}
\newcommand{\Xo}{X^\circ}
\newcommand{\tV}{\tilde{V}}

\newcommand{\fs}{f^\#}
\newcommand{\tf}{\tilde{f}}
\newcommand{\Qf}{\mathfrak{Q}}
\newcommand{\zf}{\mathfrak{z}}
\newcommand{\tfs}{\tilde{f}^\#}
\newcommand{\Tc}{\mathcal{T}}
\newcommand{\Uc}{\mathcal{U}}
\newcommand{\PGL}{\mathrm{PGL}}
\newcommand{\PK}{\mathrm{PK}}
\newcommand{\Syme}{\mathrm{Sym}_2}
\newcommand{\rr}{d}
\newcommand{\ep}{\epsilon}
\newcommand{\Qo}{Q^\circ}
\allowdisplaybreaks[4]

\begin{document}
\title{Hilbert Eisenstein Series as Doi-Naganuma Lift}
\author[Y.~Li]{Yingkun Li}
\author[M.~Zhang]{Mingkuan Zhang}
 \address{
Max Planck Institute for Mathematics,
    Vivatsgasse 7, 
    D--53111     Bonn,
    Germany}
\email{yingkun@mpim-bonn.mpg.de}
\address{Fachbereich Mathematik,
Technische Universit\"at Darmstadt, Schlossgartenstrasse 7, D--64289
Darmstadt, Germany}
\email{mzhang@mathematik.tu-darmstadt.de}
\subjclass[2020]{}
\thanks{
}

\begin{abstract}
  In this paper, we show that incoherent Hilbert Eisenstein series for a real quadratic fields can be expressed as the Doi-Naganums lift of an incoherent Eisenstein series over $\Qb$. As an application, we show when $N$ is odd and square-free, the values at Heegner points of Borcherds product on $X_0(N)^2$ with effective divisors are not integral units when the discriminants are sufficiently large.
  This generalizes the main result in \cite{Li21} to higher levels.
  In the process, we explicitly describe the Rankin-Selberg type $L$-function that appeared in the work of Bruinier-Kudla-Yang \cite{BKY12} when the quadratic space has signature $(2, 2)$, and give a new construction of fundamental invariant vectors appearing in Weil representations of finite quadratic modules. 
\end{abstract}
\date{\today}
\maketitle

 \makeatletter
 \providecommand\@dotsep{5}
 \def\listtodoname{List of Todos}
 \def\listoftodos{\@starttoc{tdo}\listtodoname}
 \makeatother

\tableofcontents

\section{Introduction}
\label{sec:intro}
Let $F$ be a real quadratic field and  $\chi$ a finite order Hecke character of $F$.
Denote $I(s, \chi)$ the set of induced sections for $G_F := \mathrm{Res}_{F/\Qb}G$, where $G = \SL_2$.
Given a standard section $\Phi(h, s) \in I(s, \chi)$, one can form the Hilbert Eisenstein series
\begin{equation}
  \label{eq:Eh}
  E(h, \Phi, s) := \sum_{\gamma \in B(F)\backslash G(F)} \Phi(\gamma h, s),
\end{equation}
which is an automorphic form on the Hilbert modular surface associated with $F$.
For suitable choices of $\Phi$, one recovers the  classical holomorphic Eisenstein series as in \cite{Freitag}.
They are crucial ingredients in the study of special values of $L$-functions,

When $\chi$ is quadratic and ramifies at all the archimedean places, it corresponds to a CM quadratic extension $\k/F$.
In this case, one can view $\k$ as an $F$-quadratic space $W_\alpha$ with respect to the quadratic form $Q_\alpha := \alpha \cdot \Nm_{\k/F}$ for $\alpha \in F^\times$,
and construct standard sections in $I(s, \chi)$ from Schwartz functions on $W_\alpha(\Ab_F)$
using the Weil representation.
%
They are called Siegel sections (see Definition \ref{def:Siegel-sec}),
and the associated Eisenstein series play an important role in the Siegel-Weil formula.
When the archimedean components of $\varphi$ are Gaussians, the corresponding Hilbert Eisenstein series will have weight $(\sgn(\alpha), \sgn(\alpha'))$, where $\prime$ denotes conjugation in $\Gal(F/\Qb)$.


Suppose from now on $\Nm(\alpha) < 0$.
Then we can alter the sign of the Hasse invariant of
$W_\alpha$ at the archimedean place where $\alpha$ is negative,
and construct a standard section $\Phi \in I(s, \chi)$ that does not come from a global $F$-quadratic space.
Such a section is called \textit{incoherent} following \cite{Kudla97}, and the corresponding  Eisenstein series $E(h, \Phi, s)$ vanishes at $s = 0$.
Its first derivative at $s = 0$ plays a crucial role in the work of Gross and Zagier on singular moduli \cite{GZ85}, the Gross-Zagier formula \cite{GZ86}, and ultimately the Kudla program \cite{Kudla97}.
Note that the incoherent Eisenstein series has parallel weight $1$.

In a classical work \cite{DN69}, Doi and Naganuma gave a construction of holomorphic Hilbert modular forms of parallel weight from elliptic modular forms when $F$ is real quadratic.
On the one hand, this is the first instance of base change \cite{Saito75, Shintani79, Langlands80}.
On the other hand, this is an example of theta lifting from $G_\Qb$ to $\mathrm{O}(2, 2)$, which is isogenous to $G_F$ \cite{Zagier75, Kudla78}.
The second perspective is particularly useful in understanding the diagonal restrictions of Hilbert modular forms.
In the literature, this lifting map mostly deals with the case for holomorphic modular forms having parallel weight $k$ with $k$ sufficiently large.
This in particular excludes the case of Hilbert Eisenstein series having small or non-parallel weight. 
One of the goals in this paper is to fill this void and demonstrate that incoherent Eisenstein series for $G_F$ with real quadratic $F$
are Doi-Naganuma lifts of Eisenstein series on $G_\Qb$ when $\k/\Qb$ is Galois.

%
To be precise, let $V \cong \Qb^2 \oplus F$ be a rational quadratic space 
and $\varphi \in \Sc(V(\Ab))$ a Schwartz function, we can associate a theta kernel $\theta_V(g, h, \varphi)$ (see section \ref{subsec:V}).
%
For a Galois extension $\k/\Qb$ of degree $4$, there exist two characters $\chi_1, \chi_2$ of $\Qb$ such that $\chi_i\circ \Nm = \chi$.
The order of $\chi_i$ is the exponent of $\Gal(\k/\Qb)$.
Given a standard section $\Xi \in I(s, \chi_1)$, we form the theta integral
\begin{equation}
  \label{eq:I}
  I(h, \Xi, \varphi, s) := \int_{G(\Qb)\backslash G(\Ab)} E(g, \Xi, s) \theta_V(g, h, \varphi) dg,
\end{equation}
which converges absolutely when $\Re(s) \gg 0$.
Our first main result is as follows.
\begin{theorem}
  \label{thm:main}
Let $\k/F$ be a CM extension with $\k/\Qb$ Galois, and denote $\chi = \chi_{\k/F}$ the associated Hecke character.
  For any Siegel section $\Phi \in I(s, \chi)$,
  function $C(s)$ real-analytic at $s = 0$ and $n \ge 1$, there exist a standard section $\Xi \in I(s, \chi_2)$ and $\varphi \in \Sc(V(\Ab))$ such that
  \begin{equation}
    \label{eq:main}
  I(h, \Xi, \varphi,s ) = C(s) E(h, \Phi, s) + O(s^n).
  \end{equation}
\end{theorem}


In specific cases of $\Phi$, we can choose suitable $C(s)$ and
be more precise about $\Xi, \varphi$ and the error term. For example, let  $D_1, D_2 < 0$ be
fundamental discriminants and $D := D_1D_2$ not a perfect square.
Then
$\k = \k_1\k_2$
is a biquadratic extension of $\Qb$ with $\k_i = \Qb(\sqrt{D_i})$, and contains the real quadratic field  $F = \Qb(\sqrt{D})$.
All of $I(s, \chi_2)$ are Siegel sections from $\Sc(\Ab_{\k_2})$, and we use $I(h, \varphi_2 \otimes \varphi_1, s)$ to denote 
$I(h, \lambda(\varphi_2), \varphi_1, s)$ for $\varphi_2 \in \Sc(\Ab_{\k_2}), \varphi_1 \in \Sc(V(\Ab))$. 
In this setting, we have the following result.
\begin{theorem}
  \label{thm:ex}
  Let $D_1, D_2 < 0$ be co-prime fundamental discriminants such that $D_1$ is odd.
For $N = k = 1$, let  $\Phi^{(k, k)}_N   =
\Phi^\infty_N
\otimes \Phi^{(k, k)}_\infty \in I(s, \chi)$ be the incoherent Siegel section associated with $\cha(L_N) \in \Sc(\widehat \Vo) \cong \Sc(\hat W)$
  (see \textsection \ref{subsec:Hilbert-Eisen} and  \textsection \ref{subsec:CM} for notations). 
  Choosing
$\tv^{k}_N = \tv_N^\infty \otimes (\varphi_\infty \otimes \varphi_\infty^{(k, 0)})$ with
$\tv_N^\infty \in \Sc(\widehat{\k_2} \oplus \hat V)$ as in \eqref{eq:tvf},
$\varphi_\infty \in \Sc(\k_2 \otimes \Rb)$ the standard Gaussian and $\varphi^{(k, 0)}_\infty \in \Sc(V(\Rb))$ the Schwartz function defined in \eqref{eq:varphikl} gives us
  \begin{equation}
    \label{eq:match-Eis}
4    I(h, \tv_N^{k}, s)
=
\Lambda(s+1, \chi_1)
E(h, \Phi^{(k, k)}_N, s),
\end{equation}
where $\Lambda(s, \chi_1)$ is the completed $L$-function associated with  $\chi_1$ (see \eqref{eq:E*}).
\end{theorem}

\begin{remark}
  For square-free $N$, odd $k \ge 1$ and fundamental discriminants $D_i$ not both even, the general version of \eqref{eq:match-Eis} is in \eqref{eq:matchN}.
\end{remark}

Such an integral expression for the Hilbert Eisenstein series
leads to a much better understanding of the spectral expansion of its diagonal restrictions, which involves generalized Rankin-Selberg $L$-series.
Such $L$-series appear in Theorem 1.1 of \cite{BKY12}, where a conjecture in the spirit of the Gross-Zagier formula relating the derivative of this $L$-series and certain arithmetic intersection was given \cite[Conjecture 5.5]{BKY12}.
This was subsequently proved in \cite{AGHMP18}, which led to a proof of the averaged Colmez conjecture.
On the other hand, this $L$-series is not a standard Rankin-Selberg integral as it involves the pull-back of the Hilbert Eisenstein series. Explicit description of this $L$-series was hoped for in \cite{BKY12}, yet has not been worked out in general.

For certain coherent Hilbert Eisenstein series of weight $(1, 1)$ over a real quadratic field, the first author was able to compute the spectral expansion in Theorem 1.5 of \cite{Li17}. The method is to relate the pull-back of the Hilbert Eisenstein series to the Shimura lift of an explicitly constructed modular form of weight $\frac32$ (see \cite{Hu16} and \cite{KP19} for alternative ways to obtain analogous spectral expansions).
Using the theta integral expression in Theorem \ref{thm:main}, we can explicitly describe this $L$-series, defined in \eqref{eq:LsG}, in the following situation.

\begin{theorem}
  \label{thm:L-func}
  Let $D_1, D_2< 0$ be distinct fundamental discriminants, not both even.
 For a newform $\Go \in S_{2k}(N)$ with odd and  square-free level $N \in \Nb$  
satisfying
  \begin{equation}
    \label{eq:Heegner}
   \kro{D_i}{p} = 1,~ \text{ for all }i = 1, 2,~ p \mid N \text{ prime factor,}
 \end{equation}
 let  $L(s, \Go; D_1, D_2)$ be the $L$-function defined in \eqref{eq:LsG}.
 Then
  \begin{equation}
    \label{eq:L-func}
    \begin{split}
          L(s, \Go; D_1, D_2)
      &=
   C_k(s) 
   \frac{ {\zeta_N(s+1)}}{ {\zeta_N(1)}}
        \frac{|\Go|_{\mathrm{Pet}}^2}{3|\tGo|_{\mathrm{Pet}}^2}
        \frac{(1 + \epsilon(\Go))\overline{c_{\tGo}(|D_1|)}{c_{\tGo}(|D_2|)}   L(\Go, s + k)
}
        {\Lambda(s+1, \chi_1)\Lambda(s+1, \chi_2 )|D_1D_2|^{(k-1)/2}
   }
    \end{split}
  \end{equation}
  where $C_k(s)$ is the constant defined in \eqref{eq:LsG0}, 
  $\zeta_N(s)$ is the partial Riemman $\zeta$-function, $\tGo = \sum_{m > 0} c_{\tGo}(m)q^m \in S^+_{k+1/2}(4N)$ is the Shintani lift of $G$, and $\epsilon(\Go)$ is the Fricke eigenvalue of $\Go$.
\end{theorem}

\begin{remark}
  \begin{enumerate}
  \item   More general result for an arbitrary  eigenform $\Go$ of higher weight is in Theorem \ref{thm:LsG}.

\item
    The sign of the functional equation of $L(\Go, s)$ at $s = k$ is $-\epsilon(\Go)$. 
    When $\epsilon(\Go) = -1$,   the $L$-function $L(s, \Go; D_1, D_2)$ vanishes identically.
    Otherwise, it vanishes at $s = k$.
    Both cases agree with the vanishing of the incoherent Eisenstein series at $s = 0$. 
  \end{enumerate}

\end{remark}

Finally, we give an arithmetic application of such a spectral expansion.
In \cite{Li21}, the first author showed that the difference of any two singular moduli is never an integral unit.
This gives a different proof and generalizes the result in \cite{BHK}, where one of the singular moduli is 0.
The proof utilizes the result of Gross-Zagier on factorization of norm of differences of singular moduli \cite{GZ85}, and generalization by Gross-Kohnen-Zagier to higher weight \cite{GKZ87}.

On $X_0(N)^2$ for a higher level modular curve $X_0(N)$ , the generalizations of $j(z_1) - j(z_2)$ are Borcherds products $\Psi_f$ associated with weakly holomorphic function $f$ on $X_0(N)$ (see section \ref{subsec:hGf}).
It is defined over $\Zb$ if $f$ has integral Fourier coefficients and a meromorphic function if $f$ has no constant term. Its divisor is determined by the principal part Fourier coefficients of $f$, which is effective precisely when these coefficients are non-negative. 
In this case, their values at Heegner points on $X_0(N)^2$ are algebraic integers.
The strategy in \cite{Li21} to prove that they are not algebraic units runs into problem due to the appearance of $L'(0, \Go; D_1, D_2)$.
With Theorem \ref{thm:L-func} and standard subconvexity results for $\GL_1$-twists of a $\GL_2$ automorphic form, we can overcome this problem and prove the following higher level generalization of the main result in \cite{Li21}.

\begin{theorem}
  \label{thm:sing}
  Let $N\in \Nb$ be odd, square-free, and
  $\Psi(z_1, z_2)$ be a Borcherds product on $X_0(N)^2$ defined over $\Zb$ with effective divisor. 
  For Heegner points $\tau_1, \tau_2$ on $X_0(N)$ with fundamental discriminants $D_1, D_2 < 0$ satisfying \eqref{eq:Heegner} and not both even, the algebraic integer $\Psi(\tau_1, \tau_2)$ is not an integral unit 
  when  $\max(|D_1|, |D_2|)$ is sufficiently large. 
\end{theorem}


\begin{remark}
  The proof of Theorem \ref{thm:sing} uses Siegel's class number lower bound, and is therefore not effective. 
\end{remark}
\begin{remark}
  The condition for $N$ to be odd and square-free is technical, and can be removed with some work. 
\end{remark}
\begin{remark}[Genus 0]
We give an example of Theorem \ref{thm:sing}. 
Let $\pi_N(\tau)$ be the hauptmodul of a genus zero modular curve $X_0(N)$ with $N$ odd and square-free, i.e.\ $N \in \{3, 5, 7, 13\}$.
Then $\pi_N(z_1) - \pi_N(z_1)$ is a Borcherds product. 
Suppose $D_1$ and $D_2$ are distinct fundamental discriminants, not both even, and coprime to $N$.
For Heegner points $\tau_i$ of discriminant $D_i$ on $X_0(N)$, the previous theorem shows that the algebraic integer
$\pi_N(\tau_1)-\pi_N(\tau_2)$ is not an integral unit when $D_1D_2$ is sufficiently large.
\end{remark}

\subsection*{Proof Strategy and Innovations}
The proof of Theorem \ref{thm:main} follows from unfolding the left hand size of the equality \eqref{eq:main} and reduce it to a local matching problem.
This is contained in section \ref{subsec:V}, where we also describe the local section from $I(h)$ as integral of a certain local invariant vector (see \eqref{eq:Fcp}). 
At the archimedean and unramified finite places, the matching results are standard and covered in sections \ref{subsec:arch} and beginning of section \ref{subsec:non-arch-I}.
For the other places, the matching data $\Xi$ and $\varphi$ in \eqref{eq:main} are not explicitly given in terms of $C(s)$ and $\Phi$. Instead, we use the decomposition of Siegel sections in \cite{KR92} and non-vanishing of Whittaker coefficients to show the existence of matching data.
This is contained in the second half of section \ref{subsec:non-arch-I} and \ref{subsec:non-arch-II}. 
Analogous argument was used in \cite{BLY22} to prove matching result when the input to the theta lift is a cusp form instead of an Eisenstein series.

To deduce the explicit example in Theorem \ref{thm:ex}, or more generally \eqref{eq:match-Eis}, we explicitly write down the Schwartz function $\varphi$ and section $\Xi$ at the ramified places, i.e.\ those dividing $DN$, in section \ref{sec:exs}. At the place dividing $D_1$, we need certain fundamental invariant vectors in the local Weil representation. These can be traced back at least to the pioneering work of Shintani \cite{Shin75}, and also appeared in the classical work by Gross-Kohnen-Zagier \cite{GKZ87}. It adds a quadratic twist to the theta kernel, which produces the special values of quadratic twist $L$-functions.
Such invariant vector plays an important role in constructing theta functions \cite{BO10, BEY21}, and have been studied extensively in the context of vector-valued modular forms \cite{ES17, Zemel23, Bie23, MS22}. 
In section \ref{subsec:inv-vec}, we give a new way to construct these fundamental invariant vectors using a quadratic character $\chi$ of the \textit{special orthogonal group} of the finite quadratic module. In fact, we show that the $\chi$-isotypic part of the group ring for the Weil representation is 1-dimensional, from which one easily deduces the $\SL_2(\Zb_p)$-invariance of these vectors. This can be viewed as the Howe correspondence on the level of finite Weil representations. 
It would be really interesting to extend this to other finite quadratic modules over valuation rings of finite extensions of $\Qb_p$.

To obtain \eqref{eq:L-func} in Theorem \ref{thm:L-func}, or more generally Theorem \ref{thm:LsG}, we substitute in the theta integral expression for the Hilbert Eisenstein series and unfold the integral defining $L(s, \Go; D_1, D_2)$. To work out the local factors at places dividing $N$, we prove an explicit version of the twisted Shintani theta lift of eigenforms, in particular oldforms.
This is the content of Proposition \ref{prop:Shintani}, and will be very helpful for future works involving explicit Shintani lifts. 
With Theorem \ref{thm:LsG} in hand, the proof of Theorem \ref{thm:sing} is a follows from adapting the proof in \cite{Li21} with the addition of equidistribution property of CM points \cite{Duke88}. 

\subsection*{Acknowledgement}
Theorem \ref{thm:ex} originated during earlier discussions between the first author and Stephan Ehlen, which we heartily thank him for. 

Both authors are supported by the Deutsche Forschungsgemeinschaft (DFG) through the Collaborative Research Centre TRR 326 ``Geometry and Arithmetic of Uniformized Structures'' (project number 444845124). 
The first author is also supported by the Heisenberg Program ``Arithmetic of real-analytic automorphic forms''(project number 539345613).

This work was completed while the first author was visiting MPIM Bonn. He thank the institute for its friendly working environment. 

\section{Preliminaries}
\label{sec:prelim}
Some notations: for $\alpha \in \Cb$, denote $\ebf(\alpha) := e^{2\pi i \alpha}$.
For a number field $K$, denote $\df_K, \Oc_K, \Ab_K$ and $\hat K$ its different, ring of integers, adeles and finite adeles respectively. For $K = \Qb$, we write $\Ab:= \Ab_\Qb$.
Let $\psi = \psi_\Qb = \otimes_p \psi_p : \Ab/\Qb \to \Cb$ be the normalized additive character such that $\psi_\infty(x) = \ebf(x)$.
For any number field $K$, we define $\psi_K: \Ab_K/K \to \Cb$ by $\psi_F(x) := \psi(\tr^K_\Qb(x))$.
For a $\Zb$-module $R$ and any $N \in \Nb$, we denote $R_N :=  R \otimes_\Zb\prod_{p \mid N \text{ prime}}\Zb_p \subset \hat R := R \otimes \hat \Zb$. 
Similarly for $a \in R$, we write $a_N := a \otimes 1 \in R_N \subset \hat R$. 
\subsection{Group and Measures}
For any commutative ring $R$, denote as usual
$$
M(R) := \{m(a) = \smat{a}{}{}{a^{-1}}: a \in R^\times\} \subset \SL_2(R),~
N(R) := \{n(b) = \smat{1}{b}{}{1}: b \in R \}\subset \SL_2(R).
$$
and for $N \in \Nb$
\begin{equation}
  \label{eq:PK}
  \begin{split}
    \PK(N) 
    &  = \prod_{p < \infty}    \PK_p(N) 
      := \{\gamma \in \PGL_2(\hat\Zb): \gamma \equiv \smat{*}{*}{0}{*} \bmod N\},~\\
    \K(N) &
= \prod_{p < \infty}    \K_p(N) 
            := \{\gamma \in \SL_2(\hat\Zb): \gamma \equiv \smat{*}{*}{0}{*} \bmod N\}.
  \end{split}
\end{equation}
If $p \nmid N$, resp.\ $N = 1$, then we omit $N$ from the notation in $\PK_p(N)$ and $\K_p(N)$, resp.\ $\PK(N)$ and $\K(N)$.


The Haar measure on $\SL_2(\Rb)$
is given by $d\mu(\tau)\frac{d\theta}{2\pi} = \frac{dudv}{v^2}\frac{d\theta}{2\pi}$ with the coordinate
\begin{equation}
  \label{eq:sl2R}
  g_\tau\kappa(\theta) \in \SL_2(\Rb),~
  g_\tau := n(u)m(\sqrt{v}),~
  \kappa(\theta) := \smat{\cos \theta }{\sin \theta }{-\sin \theta }{\cos \theta }  .
\end{equation}
The Haar measure $dg_p$ on $\SL_2(\Qb_p) = M(\Qb_p)N(\Qb_p)\K_p$
is given by
\begin{equation}
  \label{eq:measures}
  dg_p = dm(a) dn(b) dk = d^\times a db dk
\end{equation}
 with $d^\times a = \frac{da}{(1-p^{-1})|a|_p}$, $da$ and $dk$, Haar measures on $\Qb_p^\times$, $\Qb_p$ and $\K_p$, normalized such that $\Zb_p^\times, \Zb_p$ and $\K_p$ have volumes 1.
Note that $\K_p$ has the decomposition
\begin{equation}
  \label{eq:Kp}
  \K_p = J_p \sqcup_{j = 1}^p n(j)w J_p = N_p M_p (J_p \cap N^-_p) \sqcup N_p M_p N_p^- w,
\end{equation}
where $w := \smat{}{1}{-1}{}$,
$N_p := N(\Zb_p),~ M_p := M(\Zb_p),~ N^-_{p} := wN_pw^{-1}$,
and
\begin{equation}
  \label{eq:Jp}
  J_p :=  \{\smat{*}{*}{c}{*} \in \K_p: c \in p\Zb_p\} = N_p M_p (J_p \cap N_p^-)
\end{equation}
is the Iwahori subgroup.
In the coordinate $k = n(b)m(a)n^-(c) \in J_p$ with $n^-(c) := \smat{1}{}{c}{1}$, we have
\begin{equation}
  \label{eq:dk}
  dk =
  \frac{db d^\times a dc}{(1+p^{-1})|a|_p^2}
=   \frac{db d a dc}{(1-p^{-2})|a|_p^3}
  .
\end{equation}
The Haar measure on $\SL_2(\Ab)$ is the product of these local measures.

\subsection{Hecke Characters}
\label{subsec:Hecke}
Let $\k/\Qb$ be a totally imaginary, quartic Galois extension.
Then it contains a unique real quadratic field $F$, whose discriminant is denoted by $D$ and corresponds to a quadratic Hecke character $\chi_{F}$ of $\Qb$.
The quadratic Hecke character $\chi = \chi_{K/F}$ of $F$ associated with $K$ is related to a Hecke character of $\Qb$.
There are two cases to consider, either $\k$ is biquadratic or cyclic.

In the biquadratic case, let $\k_i/\Qb$ be distinct imaginary quadratic fields in $\k$ with discriminants $D_i < 0$.
Then $DD_1D_2 \in \Nb$ is a perfect square.
Denote $\chi_i := \chi_{\k_i}$ the Hecke character  corresponding to $\k_i/\Qb$.
Note that $\chi_i(a) = (a, D_i)_\Qb$ as a character on $\Ab^\times/\Qb^\times$, where  $(,)_\Qb$ is the Hilbert symbol.
In the cyclic case, let $\chi_1$ be a Hecke character of $\Qb$ of order 4 such that $\ker(\chi_1) = \k$. Note that the kernel of the character $\chi_2 := \overline{\chi_1}$ is also $\k$.

In both cases, we have
\begin{equation}
  \label{eq:char-prod}
  \chi_1 \chi_F = \chi_2,~   \chi_2 \chi_F = \chi_1,~
\end{equation}
and the following result.
\begin{lemma}
  \label{lemma:chi1}
  Let $\chi, \chi_i$ be as above. Then we have
  \begin{equation}
    \label{eq:chi1}
    \chi(a) = \chi_i(\Nm(a))
  \end{equation}
  for all $a \in \Ab_F^\times$ and $i = 1, 2$.
\end{lemma}

\begin{proof}
In the cyclic case, by viewing $\chi,\chi_i$ as characters of $\Gal(K/F),\Gal(K/\Qb)$ respectively, the restriction of $\chi_i$ to $\Gal(K/F)$ is also a quadratic character, which must equal to $\chi$. From Artin's reciprocity, we have  $\chi = \chi_i\circ\Nm$. In the biquadratic case, the restriction of $\Gal(K/F)$ to $\Gal(K_i/\Qb)$ is an isomorphism and \eqref{eq:chi1} can be proved similarly.
\end{proof}
Denote $\zeta_L(s)$ the Dedekind zeta function for a number field $L$, and $L(s, \rho)$ the $L$-function associated with $\rho \in \{\chi, \chi_1, \chi_2, \chi_F\}$.
When $\k/\Qb$ is quartic and Galois, it is easy to check that
\begin{equation}
  \label{eq:Lfac}
\zeta(s) L(s, \chi_F) L(s, \chi_1) L(s, \chi_2) =  \zeta_K(s) = \zeta_F(s) L(s, \chi) = \zeta(s) L(s, \chi_F) L(s, \chi) .
\end{equation}

\subsection{Hecke Operators on Adelic Automorphic Forms}
For a classical modular form $f \in S_{2k}(N)$ of even weight $2k$ and level $\Gamma_0(N)$, let $\langle, \rangle_{\mathrm{Pet}}$ be the Petersson inner product induced by the Petersson norm
\begin{equation}
  \label{eq:Pet-norm}
  \|f\|_{\mathrm{Pet}} := \int_{\Gamma_0(N)\backslash \Hb} |f(z) |^2 y^{2k} d\mu(z),
\end{equation}
This differs from the definition in \cite{Kohnen-FC} by the factor $[\SL_2(\Zb):\Gamma_0(N)]^{-1}$.
\footnote{For half-integral weight modular form of level $4N$, we define $  \|\cdot\|_{\mathrm{Pet}}$ in the same way without the factor $[\SL_2(\Zb):\Gamma_0(4N)]^{-1}$.}
For $m\in \Nb$ co-prime to $N$, let $T_m$ be the classical Hecke operator, which is self-adjoint with respect to  $\langle, \rangle_{\mathrm{Pet}}$ \cite[section 5.2]{DS-book}.

We use $\fs$ denote its adelization, which is an automorphic form on $\PGL_2(\Ab)$ that satisfies
 (see \cite[Prop.\ 1.4]{Kudla-book-chap})
\begin{equation}
  \label{eq:fs}
  \begin{split}
      \fs(\gamma h \kappa) &= \fs(h),~ 
    \fs(h_\infty) = f(h_\infty \cdot i) j(h_\infty, i)^{-2k},\\
    j(\smat{a}{b}{c}d, z) &:= (cz + d) (ad-bc)^{-1/2}
  \end{split}
\end{equation}
for all $\gamma \in \PGL_2(\Qb), h \in \PGL_2(\Ab), h_\infty \in \PGL^+_2(\Rb) \subset \PGL_2(\Ab)$ and $\kappa \in \PK(N)$. 
For a prime $p$ co-prime to $N$ and compactly supported, $\PK_p$-biinvariant function $\phi$ on $\PGL_2(\Qb_p)$, the associated Hecke operator $\Tc_\phi$ is defined by
\begin{equation}
  \label{eq:Tphi}
  (  \Tc_\phi \fs)(h)
  := \int_{\PGL_2(\Qb_p)} \fs(h \tilde h^{-1}) \phi(\tilde h) d \tilde h. 
\end{equation}
If $\phi = (p+1)\cha(\PK_p \smat{p}001 \PK_p)$, we denote $\Tc_\phi$ by $\Tc_p$, which is explicitly given by
\begin{equation}
  \label{eq:Tcp}
  (  \Tc_p \fs)(h) =
  \sum_{\beta \in
\PK(p) \backslash \PK}
\fs(h (\smat{p}001_p \beta)^{-1})
=
p^{1-k}(T_p f)^\#(h).
\end{equation}
Also for $d \mid N$ such that $\gcd(d, N/d) = 1$, the adelic Atkin-Lehner operator $\Wc_d$ is given by
\begin{equation}
  \label{eq:Wcd}
  \Wc_d \fs
  :=
\rho\lp \smat{}{-1}d{}_d  \rp  \fs,
\end{equation}
where we denote for any $h_0 \in \PGL_2(\Ab)$
\begin{equation}
  \label{eq:rho}
(  \rho(h_0) \fs)(h) := \fs(hh_0).
\end{equation}
Also for $N'\in \Nb$, we have the operator $(V_{N'}f)(z) := f(N'z)$, whose adelic incarnation is
\begin{equation}
  \label{eq:Vd}
  (V_{N'}f)^\# = (N')^{-k} \rho\lp \smat{1/N'}{}{}1_{N'} \rp \fs. 
\end{equation}
Now for any $N' \mid N$, we define the trace operator on adelic modular forms via
\begin{equation}
  \label{eq:Tr}
  \tr^N_{N'} \fs := \sum_{\kappa \in \PK_{}(N') / \PK_{}(N)  } \rho(\kappa) \fs,
\end{equation}
and on classical modular forms through their adelizations.
For $p\mid N$ and $\fs$ of level $N$, we can define $\Tc_p$ on $\fs$ via \eqref{eq:Tcp}. 
If $\gcd(N/p, p) = 1$, then
\begin{equation}
  \label{eq:tr-Hecke}
\tr^{N}_{N/p} \fs = \Tc_p\lp \rho\lp \smat{p}{}{}{1}_p \rp\fs \rp
\end{equation}
is a direct consequence of \eqref{eq:Tcp}.

Let $G'_{\Ab}$ be the metaplectic cover of $\SL_2(\Ab)$ \cite[Eq.\ (0.11)]{Kudla96}, and denote $G'_R$ the preimage of $\SL_2(R)$ for any subring $R \subset \Ab$. 
For a half-integral weight modular form $\tf \in S_{k+1/2}(4N)$ of level $\tilde\Gamma_0(4N) \subset G'_{\Zb}$, one can adelize it to an automorphic form $\tfs$ on $G'_{\Qb}\backslash G'_{\Ab}$ satisfying an analogue of \eqref{eq:fs} (see e.g.\ Lemma 1.1 in \cite{Kudla03}). 
For $m \in \Nb$ co-prime to $4N$, we have the following operators
\begin{equation}
  \label{eq:Um}
  (  U_m \tf )(\tau) := 
\frac1m  \sum_{j = 1}^{m} \tf \lp \frac{\tau + j}{m} \rp.
\end{equation}
If $m = \rr^2$, then $U_{\rr^2} \tf$ is modular of weight $k+1/2$ and level $4N\rr$.
It is easy to check that
\begin{equation}
  \label{eq:Ucm}
\rr^{-k-1/2}
  (U_{\rr^2} \tf)^\#(g') = (\Uc_{\rr^2} \tfs)(g')
  :=
\rr^{-2}
  \sum_{j = 1}^{\rr^2} \tfs \lp g' [(m(\rr)n(j))_{\rr}, 1] \rp.
\end{equation}
In \cite{Kohnen-newform}, Kohnen showed that
$S_{k+1/2}(4N)$ contains the plus subspace $S^+_{k+1/2}(4N)$ and subspace $S^{+, \mathrm{new}}_{k+1/2}(4N)$ spanned by newforms when $N$ is square-free. Furthermore, 
\begin{equation}
  \label{eq:Kohnen-isom}
S^+_{k+1/2}(4N) \cong
  \oplus_{\rr \mid N} U_{\rr^2}(S^{+, \mathrm{new}}_{k+1/2}(4N/\rr))
\end{equation}
as Hecke modules.

\subsection{Weil Representation}
For a rational quadratic space $(V, Q)$ of even dimension $m$, denote $\omega = \omega_\psi$ the global Weil representation of $\SL_2(\Ab) \times H(\Ab)$, with $H = H_V := \mathrm{Spin}_V$.
On the space of Schwartz functions $\Sc(V(\Ab))$, it is given by the explicit formula
\begin{align*}
  (\omega(m(a))\varphi)(x)
  &= |a|^{m/2} \chi_V(a) \varphi(ax),
  (\omega(n(b))\varphi)(x)
  = \psi(bQ(x)) \varphi(x),\\
  (\omega(w) \varphi)(x)
  &= \int_{V(\Ab)} \varphi(y) \psi((x, y)) dy,
  (\omega(h) \varphi)(x)
  = \varphi(h^{-1}x).
\end{align*}
The measure $dy$ on $\Ab$ is self-dual with respect to $\psi$.
Locally, $\SL_2(\Qb_p) \times H_V(\Qb_p)$ acts on $\Sc(V(\Qb_p))$ via the Weil representation $\omega_p$ with analogous formula.
For $\varphi \in \Sc(V(\Qb_p))$, denote
\begin{equation}
  \label{eq:Kp-inv}
  \varphi^{\K_p} := \int_{\K_p} \omega_p(k) \varphi dk
\end{equation}
the $\K_p$-invariant component of $\varphi$.

Though $\Sc(V(\Qb_p))$ is infinite dimensional, the action of $\omega$ restrict to the subspace
\begin{equation}
  \label{eq:ScL}
  \Sc_L := \{\varphi \in \Sc(\hat V): \mathrm{supp}(\varphi) \subset \hat L^\vee, \varphi(x + \lambda) = \varphi(x) \text{ for all }x \in \hat V, \lambda \in \hat L\}
\end{equation}
for any even, integral lattice $L \subset V$ with dual lattice $L^\vee$ and completion $\hat L := L \otimes \hat\Zb$.
The natural inclusion $L \subset \hat L$ induces canonical isomorphism $L^\vee/L \cong \hat L^\vee/\hat L$, and we identifies them in this way.
The subspace $\Sc_L$ is finite dimensional and isomorphic to $\Cb[L^\vee/L]^\vee$.
Let  $\langle, \rangle$ be the $\Cb$-bilinear pairing on $\Cb[L^\vee/L]$ given by
\begin{equation}
  \label{eq:C-pair}
  \langle\ef_\nu, \ef_\mu \rangle
  :=
  \begin{cases}
    1& \text{ if } \nu = \mu,\\
    0 &\text{ otherwise,}
  \end{cases}
\end{equation}
where $\{\ef_\mu\}_{\mu \in \hat L^\vee/\hat L}$ is a basis of $\Cb[\hat L^\vee/\hat L]$.
Let $\rho_L = \rho_{L^\vee/L}$ denote the action of $\SL_2(\Zb) \subset\SL_2(\hat\Zb)$ on $\Cb[L^\vee/L] \cong \Cb[\hat L^\vee/\hat L]$ via $\omega$ and this duality.
Then it is the Weil representation associated with the finite quadratic module $L^\vee/L$ as in \cite{Borcherds98}.
Similarly, the local Weil representation $\omega_p$ induces the representation $\rho_{L_p}$ of $\K_p$ on $\Cb[L_p^\vee/L_p]$, with $L^{\vee}_p := L^{\vee}\otimes \Zb_p$ and $L^\vee_p/L_p \cong (L^\vee/L)_p := (L^\vee/L) \otimes \Zb_p$.

\subsection{Hilbert Eisenstein Series}
\label{subsec:Hilbert-Eisen}
Let $F$ be a totally real field of degree $d$ with different $\df_F$.
For a Hecke character $\chi$ of $F$, denote $I(s, \chi)$ the functions $\Phi(h, s)$ on $G_F(\Ab)$ satisfying
$$
\Phi(n(\beta) m(\alpha)h, s) = \chi(\alpha) |\alpha|^{s+1} \Phi(h, s),~ s \in \Cb, \alpha \in \Ab^\times_F, \beta \in \Ab_F.
$$
A section $\Phi \in I(s, \chi)$ is called standard if there is an open compact $K \subset G(\hat F)$ such that for every $g \in K$, the function $\Phi(g, s)$ is independent of $s$.
It is factorizable if $\Phi = \otimes_{v \le \infty} \Phi_v$, 
where $\Phi_v \in I(s, \chi_v)$ with $\chi = \otimes_v \chi_v$ and $I(s, \chi_v)$ the local version of $I(s, \chi)$.
Given a (linear combination of) standard, factorizable section $\Phi \in I(s, \chi)$, we can form the Hilbert Eisenstein series $E(h, \Phi, s)$ as in \eqref{eq:Eh}.
For an archimedean place $v$ of $F$, $k \in \Zb$ and $\chi_v: x\mapsto \sgn(x)^k$ a character of $F_v^\times \cong \Rb^\times$, we have the standard section  $\Phi^k_\infty \in I(s, \chi_v)$ of weight $k$, satisfying
\begin{equation}
  \label{eq:keq}
\Phi^k_\infty(g\kappa(\theta) , s) = \ebf(k \theta) \Phi^k_\infty(g, s)
\end{equation}
for all $\theta \in \Rb, s \in \Cb$ and $g \in \SL_2(\Rb) \cong  G(F_v)$, and normalized such that $\Phi^k_\infty(1) = 1$.


When $\chi$ is quadratic, let $\k/F$ be the corresponding quadratic extension with relative discriminant $d_{K/F}$.
If $\chi$ is ramified at all archimedean places, then $\k$ is a CM field.
For $\alpha \in F^\times$, we have the quadratic space $W_\alpha$ and the map $\lambda_\alpha$ 
\begin{equation}
  \label{eq:SW-sec}
  \lambda_\alpha: \Sc(W_\alpha) \to I(0, \chi),~ \lambda_\alpha(\varphi)(g)
  = (\omega_\alpha(g)\varphi)(0),
\end{equation}
where $\omega_\alpha$ is the Weil representation associated with $W_\alpha$.

More generally for each idele $\alb = (\alpha_v)_{v \le \infty} \in \Ab_F^\times$, we can define an $\Ab_F$ quadratic space
$$
\Wb_{\alb} = (W_{\alpha_v})_{v \le \infty},~ W_{\alpha_v, v} = (K_v, \alpha_v\cdot \Nm_{K_v/F_v}).
$$
It is simply $(W_{\alpha} \otimes F_v)_v$ if $\alb \in F^\times$. More generally, we call $\Wb_\alb$ \textit{coherent} if it is isomorphic to such a quadratic space. We call it \textit{incoherent} if it becomes coherent after changing the Hasse invariant at a place $v$. We say these two $\Ab_F$-quadratic spaces are \textit{$v$-neighbors} of each other. 
For each place $v$ of $F$,
we have local analogue $\lambda_{\alpha_v, v}$ of the section map in \eqref{eq:SW-sec}.
Its images are denoted by $R(W_{{\alpha_v}, v})$ and called \textit{Siegel sections}.
The local sections $I(0, \chi_v)$ can be expressed explicitly in terms of Siegel sections as follows \cite[Theorem 2.1]{KR92}.
\begin{proposition}
  \label{prop:local-sec}
For any non-archimedean place $v$ of $F$, we have
\begin{equation}
  \label{eq:local-sec}
  I(0, \chi_v) = \bigoplus_{\alpha \in F_v^\times / \Nm(K_v^\times)}
R(W_{\alpha, v}).
\end{equation}
\end{proposition}

\begin{definition}
  \label{def:Siegel-sec}
Let $\chi$ be a quadratic character.  We call $\Phi \in I(s, \chi)$  a \textit{Siegel section}
if it is a linear combination of standard sections $\otimes_v \Phi_v \in I(s, \chi)$ with $\Phi_v \in I(0,\chi_v)$ Siegel sections for $v < \infty$ and $\Phi_v \in I(s, \chi_v)$ standard sections of weight $k \in \Zb$ for $v \mid \infty$.
A Siegel section $\Phi = \otimes_v \Phi_v$ is called \textit{coherent}, resp.\ \textit{incoherent}, if $\Phi_v = \lambda_{\alpha_v}(\varphi_v)$ for $\varphi_v \in \Sc(W_{\alpha_v})$ with $(W_{\alpha_v})_{v \le \infty}$ coherent, resp.\ incoherent.
\end{definition}
For a Siegel section $\Phi \in I(s, \chi)$ as above, we normalize
$E(h, \Phi, s)$ by \cite[section 4]{BKY12}
\begin{equation}
  \label{eq:E*}
  E^*(h, \Phi, s) := \Lambda(s+1, \chi) E(h, \Phi, s),
\end{equation}
where $\Lambda(s, \chi) := A^{s/2} \Gamma_\Rb(s+1) L(s, \chi)$ is the completed Dirichlet $L$-function associated with $\chi$ with $A := \Nm(\df_F d_{K/F})$ and
$\Gamma_\Rb(s) := \Gamma(s/2) \pi^{-s/2}$  the Gamma factor at infinity.
Note that
\begin{equation}
  \label{eq:Gamma-rat}
\frac{  \Gamma_\Rb(s+2r)}{  \Gamma_\Rb(s)} = \frac{s(s+2)\cdots (s+2(r-1))}{(2\pi)^r}
\end{equation}
for any $r \in \Zb_{\ge 0}$.

If $\Phi = \lambda(\cha(\hat P)) \otimes_{1 \le j \le d} \Phi_{\infty_j}^{k_j}$ is a(n in)coherent Siegel section from the ($\infty_j$-neighbor of) $F$-quadratic space $W_\alpha$ for some lattice $P \subset W_\alpha$ and  $k = (k_j)_{1\le j \le d} \in \Zb^d$, then the function
\begin{equation}
  \label{eq:EP}
  E_{P, k}(z, s) := E(h_z, \Phi, s) \prod_{1 \le j \le d} y_j^{-k_j/2}
\end{equation}
is a classical, real-analytic Hilbert Eisenstein series in $z = (z_j)_{1\le j \le d} = (x_j + i y_j)_{1 \le j \le d} \in \Hb^d$ of weight $(k_j)_{1 \le j \le d} \in \Zb^d$.
Here $h_z = (g_{z_j})_{1 \le j \le d} \in  G_F(\Rb) \subset G_F(\Ab)$. 

\subsection{Theta Function and Doi-Naganuma Lift}
\label{subsec:V}
From now on, let $F$ be a real quadratic field with discriminant $D > 0$ and conjugation $'$.
Consider the quadratic space
\begin{equation}
  \label{eq:V}
  V := \{A \in M_2(F): A^t = A'\},~ Q = \aa\cdot \det
\end{equation}
for $\aa \in \Qb_{> 0}$.
It has signature $(2, 2)$
and $\chi_V(x) = (x, D)_\Qb =  \chi_F(x)$.
For a Schwartz function $\varphi \in \Sc(V(\Ab))$, we have the theta function
$$
\theta_{}(g, h, \varphi) := \sum_{x \in V} (\omega(g) \varphi)(h^{-1} x),
$$
where $\omega = \omega_V$ is the Weil representation.
We are interested in the integral $I(h, \Xi, \varphi, s)$ defined in \eqref{eq:I} for $\Xi \in I(s, \chi_2)$, which can be expressed as an Eisenstein series on $G_F \cong H$.
\begin{lemma}
  \label{lemma:unfold}
  Let $\chi, \chi_1$ be as in Lemma \ref{lemma:chi1}.
For any $h\in H(\mathbb{A})$, we have
$$I(h,\Xi, \varphi,s)=\sum_{\gamma\in B(F)\backslash G(F)} \Fr_{}(\gamma h, s; \Xi, \varphi),$$
where $\Fc(h, s; \Xi, \varphi) \in I(s, \chi)$ is defined by the convergent integral
\begin{equation}
  \label{eq:sec}
  \begin{split}
    \Fr(h,s ) &=
    \Fr_{}(h, s; \Xi, \varphi)\\
    &:=
      \int_{(B(\hat{\Zb})B(\Rb))\backslash \SL_2(\hat{\Zb})\SL_2(\Rb)}\Xi(k)
  \int_{\Ab^\times} |a|^{s + 1}
\chi_{1}(a)
 (\omega(k) \varphi)(h^{-1}  x_0 a)  d^\times a dk
  \end{split}
\end{equation}
for $\Re(s) \ge 0$, with
 $x_0 := \smat{1}{0}{0}0 \in V$.
\end{lemma}
\begin{remark}
  \label{rmk:1-2}
  \begin{enumerate}
  \item   We can also use a section $\Xi$ in $I(s, \chi_1)$ instead. Then the character $\chi_1$ in the expression of $\Fc(h, s)$ will be replaced by $\chi_2$.
  \item If $\chi_2$ is quadratic, then the Siegel-Weil section map embeds $\Sc((\k_2 \oplus V)(\Ab)) \cong \Sc(\k_2 \otimes \Ab) \otimes \Sc(V(\Ab))$ into $I(s, \chi_2) \otimes  \Sc(V(\Ab))$.
    For $\tv \in \Sc((\k_2 \oplus V)(\Ab))$ with image $\sum_i
    \Xi_i \otimes \varphi_i \in I(s, \chi_2) \otimes \Sc(V(\Ab))$, we use $\Fc(h, s; \tv)$ to denote the section $\sum_i \Fc(h, s; \Xi_i, \varphi_i)$. 
  \end{enumerate}
\end{remark}
\begin{proof}
  First, notice that $\tilde \varphi(x) := \int_{(B(\hat{\Zb})B(\Rb))\backslash \SL_2(\hat{\Zb})\SL_2(\Rb)}\Xi(k)
  (\omega(k) \varphi)(h^{-1}  x)  dk$
  integrates over a compact region and again gives a Schwartz function on $V(\Ab)$ for any fixed $\varphi, \Xi$ and $h$.
Suppose $\tilde\varphi = \tilde\varphi_f \tilde \varphi_\infty$.
To show the convergence of
        $$
\Fc(h, s; \Xi, \varphi)=
  \int_{\Ab^\times} |a|^{s + 1}
\chi_{1}(a) \tilde\varphi(a x_0)  d^\times a,
        $$
it suffices to
prove the convergence of the integral
over $C \cap \Ab^\times$ for a compact subset $C \subset \Ab$, e.g.\ $C = N^{-1} \hat\Zb \times [0, 1]$ for fixed $N \in \Qb_{>0}$.
We can choose $N$ such that the  $\varphi_f(ax_0)$ is constant for $a \in N^{-1}\hat\Zb$.
Then
\begin{align*}
\int_{C \cap \Ab^\times} |a|^{s+1}\chi_{1}(a) \tilde\varphi(ax_0) d^\times a
&= \int_{N^{-1}\hat\Zb - \{0\}} |a|^{s+1}\chi_{1, f}(a) d^\times a
\int_0^1 a^{s+1} \tilde\varphi_\infty(ax_0) \frac{da}a,
  \end{align*}
  and
  \begin{align*}
\int_{N^{-1}\hat\Zb - \{0\}} |a|^{s+1}\chi_{1, f}(a) d^\times a
      &=
  \prod_p \int_{N^{-1}\Zb_p - \{0\}} \chi_{1, p}(a) |a|_p^{s+1} {d^\times a}
=
  \frac{N^{s+1}}{s+1} L(1+s, \chi_1)
  \end{align*}
  converges for $\Re(s ) \ge 0$ since $\chi_1$ is a quadratic character.

Using $G(\hat\Qb) = B(\Ab) \SL_2(\hat{\Zb})$, we can unfold the integral defining $I(h,\Xi,\varphi, s)$ as
\begin{align*}
  I(h, \Xi, \varphi, s)
  & = \int_{G(\Qb)\backslash G(\Ab)}
\sum_{\gamma \in B(\Qb)\backslash G(\Qb)} \Xi(\gamma g, s)
    \theta (g, h, \varphi) dg
   = \int_{B(\Qb)\backslash G(\Ab)}
    \Xi( g, s)
    \theta (g, h, \varphi) dg\\
&   = \int_{B(\Ab)\backslash G(\Ab)}
 \int_{\Qb^\times\backslash\Ab^\times}\int_{\Qb\backslash\Ab}\Xi(n(b)  m(a) g, s)\theta(n(b) m(a)g, h, \varphi)\frac{db}{|a|^2} d^\times a dg\\
&   =  \int_{(B(\hat{\Zb})B(\Rb))\backslash \SL_2(\hat{\Zb})\SL_2(\Rb)}\Xi(k)
 \int_{\Qb^\times\backslash\Ab^\times} |a|^{s + 1} \chi_2(a)\\
&\quad \times \int_{\Qb\backslash\Ab}
  \sum_{x \in V}\psi(bQ(x))\chi_{V}(a) (\omega(k) \varphi)(h^{-1} x a) db d^\times a dk.
\end{align*}
The integral on the last line vanishes unless $Q(x) = 0$, i.e.\ $x \in V$ is isotropic. The set of isotropic vectors in $V$ can be parametrized as
$$
\{\gamma^{-1} x_0 \alpha: \alpha \in \Qb^\times, \gamma \in P(\Qb)\backslash H(\Qb)\},
$$
where 
$x_0 \in V$ is the isotropic vector in Lemma \ref{lemma:unfold} and $P \subset H$ the parabolic subgroup fixing $x_0$.
Using this and \eqref{eq:char-prod}, we can furthermore write
\begin{align*}
  I(h, \Xi, \varphi, s)  &   =
\sum_{\gamma \in P(\Qb)\backslash H(\Qb)}
\int_{(B(\hat{\Zb})B(\Rb))\backslash \SL_2(\hat{\Zb})\SL_2(\Rb)}\Xi(k)
 \int_{\Ab^\times} |a|^{s + 1} \chi_1(a)
 (\omega(k) \varphi)(h^{-1} \gamma^{-1} x_0 a)  d^\times a dk\\
  &   =
\sum_{\gamma \in P(\Qb)\backslash H(\Qb)}
    \Fr_{}(\gamma h, s; \Xi, \varphi)
\end{align*}
Using $H(\mathbb{Q}) \cong G(F)$, we obtain $P(\mathbb{Q})\backslash H(\mathbb{Q}) \cong B(F)\backslash G(F)$. Hence the first statement holds. Using $n(\beta)x_0 = x_0$, $m(\alpha)^{-1} \cdot x_0  = \Nm(\alpha)^{-1} x_0$ and Lemma \ref{lemma:chi1}, it is easy to check that for $\beta \in \Ab_F, \alpha \in \Ab_F^\times$
\begin{align*}
  &\Fr(n(\beta)m(\alpha)h, s)\\
  &=  \int_{(B(\hat{\Zb})B(\Rb))\backslash \SL_2(\hat{\Zb})\SL_2(\Rb)} \Xi(k)
 \int_{\Ab^\times} |a|^{s + 1}
\chi_{1}(a)
 (\omega(k) \varphi)(h^{-1}  x_0 \Nm(\alpha)^{-1} a)  d^\times a dk\\
  &=
|\Nm(\alpha)|^{1+s}  \chi_1(\Nm(\alpha))
\int_{(B(\hat{\Zb})B(\Rb))\backslash \SL_2(\hat{\Zb})\SL_2(\Rb)}\Xi(k)
 \int_{\Ab^\times} |a|^{s + 1} \chi_1(a)
 (\omega(k) \varphi)(h^{-1}  x_0  a)  d^\times a dk\\
&=|\Nm(\alpha)|^{1+s} \chi(\alpha) \Fr(h, s)=|\alpha|_F^{1+s} \chi(\alpha) \Fr(h, s).
\end{align*}
Therefore $\Fc(h, s) \in I(s, \chi)$.
\end{proof}

Suppose $\Xi = \prod_{p \le \infty} \Xi_p$ and $\varphi = \otimes \varphi_p$ are factorizable. Then we can write $\Fr = \prod_{p \le \infty} \Fr_p$ with $\Fr_p \in \otimes_{v \mid p} I(s, \chi_v)$ defined by
\begin{equation}
  \label{eq:Frp}
  \begin{split}
      \Fr_p(h, s)
&  :=     \int_{B(\Zb_p)\backslash \SL_2(\Zb_p)}\Xi_p(k)
 \int_{\Qb_p^\times} |a|_p^{s + 1} \chi_{1, p}(a)
 (\omega_p(k) \varphi_p)(h^{-1}  x_0  a)  d^\times a dk,\\
      \Fr_\infty(h, s)
&  :=     \int_{0}^\pi \Xi_\infty(\kappa(\theta))
  \int_{0}^\infty a^{s + 1}
 (\omega_\infty(\kappa(\theta)) \varphi_\infty)(h^{-1}  x_0  a)  d^\times a \frac{d\theta}{2\pi}
  \end{split}
\end{equation}
for $\Re(s) > -1$. When $s = 0$, we omit it from the notation.

Now for $p < \infty$ and a local section $\Xi_p \in I(0, \chi_{1, p})$, we define a map $\Cc = \Cc_{\Xi_p}$ by
\begin{equation}
  \label{eq:Cc}
  \Cc:\Sc (V(\Qb_p))\rightarrow \Sc(V(\Qb_p)):
  \varphi\mapsto 
  \int_{\K_p} \Xi_p(k)(\omega_p(k)\varphi)(\cdot)dk.
\end{equation}
To see that $\Cc(\varphi)$ is Schwartz function, notice that $\Xi_p$ and the action of $\omega_p$ on $\varphi$ is locally constant.
So the integral in \eqref{eq:Cc} is a finite sum of Schwartz functions of the form $\omega_p(k)\varphi$.
In addition, the map $\Cc$ is equivariant with respect to the action of $H(\Qb_p)$.
Furthermore, when $\Xi_p = \lambda(\phi)$ is a Siegel section,
we can write
\begin{equation}
  \label{eq:Cc'}
  \Cc(\varphi)(x) = (\phi\otimes\varphi)^{\K_p}(0, x),
\end{equation}
with ${\cdot}^{\K_p}$ the invariant operator defined in \eqref{eq:Kp-inv}.

We further define
\begin{equation}
  \label{eq:Lambda}
  \Lambda_s: \Sc (V(\Qb_p)) \rightarrow I(s,\chi_p),~ \varphi\mapsto
  \lp h \mapsto \int_{\Qb_p^\times}|a|^{s+1} \chi_{1,p}(a)\varphi(ah^{-1}x_0)d^\times a\rp,
\end{equation}
which satisfies
\footnote{Here, we have used $\chi_V = \chi_F = \chi_1 \chi_2$ in \eqref{eq:char-prod}.}
\begin{equation}
  \label{eq:r-act}
  \Lambda_s(\omega_p(n)\varphi) =   \Lambda_s(\varphi),~
  \Lambda_s(\omega_p(m(r))\varphi) =
  |r|^{-s-1}\chi_{2, p}(r) \Lambda_s(\varphi)
\end{equation}
for all $n \in N(\Qb_p) \subset G(\Qb_p), r \in \Qb_p^\times$ and $h \in H(\Qb_p)$.
Then we can write 
\begin{equation}
  \label{eq:Fcp}
  \Fc_p(h, s; \Xi, \varphi) =
(\Lambda_s \circ \
  \Cc_{\Xi_p})(\varphi)(h) \in I(s, \chi_{p}) = \otimes_{v \mid p} I(s, \chi_v).
\end{equation}
Note that after multiplying by certain function in $s$, the section $\Lambda_s(\varphi)$ could become  standard. However, this is not true in general. For example, see Lemma \ref{lemma:standard} below.

When $\Xi_p$ is a Siegel section from another quadratic space $V_{2, p}$, let $\tilde V_p := V_{2, p} \oplus V_p$. 
Then for $\tv \in \Sc(\tilde V_p)$ we can write
\begin{equation}
  \label{eq:tFcp}
  \Fc_p(h, s; \tv)
  = \tilde\Lambda_s(\tv)(h),
\end{equation}
where $\tilde\Lambda_s$ is defined (analogously as in \eqref{eq:Lambda}) by
\begin{equation}
  \label{eq:tLambda}
    \tilde\Lambda_s: \Sc (\tilde V_p) \rightarrow I(s,\chi_p),~ \varphi\mapsto
  \lp h \mapsto \int_{\Qb_p^\times}|a|^{s+1} \chi_{1,p}(a)\tv^{\K_p}(0, ah^{-1}x_0)d^\times a\rp. 
\end{equation}
Here $\tv^{\K_p}$
is the $\K_p$-invariant part of $\tv$ defined in \eqref{eq:Kp-inv}. 
In particular, when $\tv = \phi \otimes \varphi$ with $\phi \in \Sc( V_{2, p})$, we have
\begin{equation}
  \label{eq:rmk-tv}
   (\Lambda_s \circ \  \Cc_{\Xi_p})(\varphi)=\tilde{\Lambda}_s(\phi\otimes\varphi),~ \Xi_p = \lambda_p(\phi).
 \end{equation}

\section{CM-Value Formula for Higher Green function}

In this section, we recall the higher Green function on $X_0(N)^2$, and the formula for its value at big CM points. The results are essentially from the classical work of Gross, Kohnen and Zagier \cite{GKZ87}.
Here, we follow the more modern approach using regularized theta lifts (see \cite{BEY21}, \cite{Li21}, \cite{BLY22}). 

\subsection{Higher Green Function}
\label{subsec:hGf}
Let $$Q_{s-1}(t):=\int_0^{\infty}\left(t+\sqrt{t^2-1} \cosh v\right)^{-s} d v$$
be the Legendre function of the second kind, which satisfies the ordinary differential equation
$$\left(1-t^2\right)\left(\partial_t\right)^2 F(t)-2 t \partial_t F(t)+s(s-1) F(t)=0 .$$
Define a function $g_s$ on $\mathbb{H}^2$ by
$$g_s(z_1,z_2):=-2 Q_{s-1}(\cosh d (z_1,z_2))=-2 Q_{s-1}\left(1+\frac{\left|z_1-z_2\right|^2}{2 y_1 y_2}\right).$$
By averaging over the $\Gamma_0(N)$-translates of the second variable, we obtain a function
$$G^N_s (z_1,z_2):= \sum_{\gamma\in\Gamma_0(N)} g_s(z_1, \gamma z_2).$$

Given any harmonic Maass form $f \in H_{2-2 k}(N)$ having pole only at $\infty$ with $k > 1$,
we can define the associated higher Green's function
by
$$G_{k, f}(z_1,z_2):=\sum_{m \geq 1} c_f(-m) m^{k-1} T_m G_k(z_1,z_2) $$
with $c_f(m)$ the $m$-th Fourier coefficient of the holomorphic part of $f$, and
$T_m$ the $m^{\text {th }}$ Hecke operator acting on $z_2$ (the same if the action is on $z_1$). The singularity of $G_{k, f}$ is
$$T_f:=\bigcup_{m \geq 1, c_f(-m) \neq 0}\left(z, T_m z\right).$$
This is the same (up to sign) as the higher Green function $\Phi^{k-1}_{\mathrm{vv}(f)}$ 
obtained from regularized theta lift of $\mathrm{vv}(f)$, where $\mathrm{vv}$ is the vector-valued lifting map in (2.47) of \cite{BLY22}
from $H_{2-2k}(N)$ to vector-valued harmonic Maass form in $H_{2-2k, \rho_{L_N}}$ for the Weil representation $\rho_{L_N}$ associated with the lattice $L_N$ in \eqref{eq:LN} (see Corollary 2.4 in \cite{BLY22}).

When $k = 1$ and $f$ is weakly holomorphic, i.e.\ $\xi(f) = 0$ with $\xi := 2iv^2 \overline{\partial_{\overline \tau}}$, the theta lift $\Phi_{\vv(f)}$ is just the Borcherds lift in \cite{Borcherds98}, which is invariant with respect to $\mathrm{O}(L_N)$.
Since $\vv(1)$ is non-trivial, we can replace $f$ by $f - c$ for a unique constant $c$ such that $\vv(f)$ has trivial constant term at the trivial coset.
Note that the Borcherds lift of $\vv(1)$ is the logarithm of an eta quotient on $X_0(N)^2$. 
Then the Borcherds lift can be written as
\begin{equation}
  \label{eq:Psif}
\Phi_{\vv(f - c)}(z_1, z_2) = 
  \log|\Psi_f(z_1, z_2)| , 
\end{equation}
where $\Psi_f(z_1, z_2)$ is a meromorphic function on $X_0(N)^2$ with divisor supported on $T_f$.
Furthermore, it has a product expansion and is called the \textit{Borcherds product} on $X_0(N)^2$ associated with $f$. 
Since $L_N$ splits off a hyperbolic plane over $\Zb$, we can inspect the Fourier expansion of $\Psi_f$ at the corresponding cusp in \cite{Borcherds98} to see that some power of it is defined over $\Zb$ when $f$ has rational Fourier coefficients at $\infty$.

For a smooth function $\phi: \mathbb{H}^2 \rightarrow \mathbb{C}$ and $k_1, k_2 \in \frac12\Zb$,
define the $r$-th Cohen operator as (see \cite[section 2.1]{BLY22})
\begin{equation}
  \label{eq:Cr}
  \begin{split}
&    \mathcal{C}_r(\phi)(\tau)
   = \frac{1}{(2 \pi i)^r} \sum_{s=0}^r(-1)^s\binom{k_1+r-1}{s} \binom{ k_2+r-1}{r-s}\left(\frac{\partial^{r-s}}{\partial \tau_1^{r-s}} \frac{\partial^s}{\partial \tau_2^s} \phi\right)(\tau, \tau)
  \end{split}
\end{equation}
When $r \ge 1$, it is clear that $\Cc_r(\phi)$ has trivial constant term. 
When $k_1 = k_2 = 1$, we have
\begin{equation}
  \label{eq:Cr11}
    \mathcal{C}_r(\phi)(\tau)   = (-4\pi)^{-r}
 (-1)^{r_0} \sum_{s = 0}^{r_0} \binom{r_0 - r - 1/2}{r_0 - s} \binom{r - r_0 - 1/2}{s} (R_{\tau_1} - R_{\tau_2})^{r-2s}(R_{\tau_1} + R_{\tau_2})^{2s} \phi \mid_{\tau_1 = \tau_2 = \tau},
\end{equation}
where $R_\tau = R_{\tau, k} = 2i \partial_\tau + \frac{k}{v}$ is the weight $k$ raising operator. 
If $f$ is modular of weight $(k_1, k_2)$, then $\mathcal{C}_r(\phi)$ is modular of weight $k_1+k_2+2 r$ (see e.g.\ section 2.2 in \cite{Li23}).

\subsection{Big CM cycles and CM-Value Formula}
\label{subsec:CM}
We now recall the big CM cycle on product of two level $N$ modular curves following \cite{YY19} and \cite{Ye22}. 

Let $\Vo$ be the quadratic space $M_2(\mathbb{Q)}$ with quadratic form $\Qo = \det$. The group $\Ho := \GSpin(\Vo)$ is a subgroup of $\GL_2^2$. For $N \in \Nb$,  we have the following lattice in $\Vo$, 
\begin{equation}
  \label{eq:LN}
  L_N := \left\{\smat{a}{b}{c}{d} \in M_2(\Zb): N \mid c\right\} \subset \Vo. 
\end{equation}
which is stabilized by $K_0(N) := \Ho(\hat\Qb) \cap \smat{\hat\Zb}{\hat\Zb}{N\hat\Zb}{\hat\Zb}^2$.
{Note that $\lp\smat{r_1}**{r_1^{-1}},\smat{r_2}**{r_2^{-1}}\rp  \in K_0(N)^2$ acts on $L_N^\vee/L_N \cong (\Zb/N\Zb)^2$ by sending $\smat{a}**d$ to $\smat{a r_1/r_2}**{d r_2/r_2}$.}
The associated Shimura variety
\begin{equation}
  \label{eq:XoN}
  \Xo_N := \Ho(\Qb) \backslash (\Hb^2 \sqcup (\Hb^- )^2) \times \Ho(\hat\Qb)/K_0(N)
\end{equation}
 is isomorphic to $Y_0(N)^2$ (see section 2.5 in \cite{BLY22}).

 Let $D_1, D_2 < 0$ be distinct fundamental discriminants, one of which is odd, and satisfying condition \eqref{eq:Heegner} with $N$ odd and square-free.
 Denote $\k_i, K, F, D$ the same as in the introduction, $\Oc_i := \Oc_{\k_i}$ for $i = 1, 2$ and $\Oc_{\k, 0} := \Oc_1\Oc_2 \subset \Oc_\k$ the subring of index $D_0 := \gcd(D_1, D_2)$. 
 Let $\tau_i = \frac{B_i+\sqrt{D_i}}{2NA_i} \in \Hb$ be  CM points such that $N\Nm(A_1 \tau_1), N\Nm(A_1 \tau_1), 2DN$ are pairwise co-prime natural numbers. 
 Then
 the embeddings $\iota_i = \iota_{A_i \tau_i}: \Oc_{i} \hookrightarrow M_2(\Zb)$ defined by 
 \begin{equation}
   \label{eq:iota-z}
   \iota_{\tau}(r) \binom{\tau}1 = \binom{r\tau}r,~
   \tau, r \in \k_i
 \end{equation}
 have image in $L_N$, are optimal and
 gives $L_N$ the structure of an $\Oc_{i}$-module.
 Note that $\iota_{\overline\tau}(r) = \iota_{\tau}(\overline r)$ and
 $(-1, \tau) \iota_\tau(\bar r) = (-r, r\tau)$.


 Let $W$ be the $F$-quadratic space $K$ with $Q_F(z)=\frac{z \bar{z}}{\sqrt{D}}$. 
We can identify $\mathrm{Res}_{F / \mathbb{Q}} W$ with $ \Vo$ via
\begin{equation}
  \label{eq:WVo}
  \phi: K \to \Vo,~ 
  \sum_{i=1}^4 x_i e_i \mapsto
\frac1N
  \left(\begin{array}{ll}
x_3 & x_1 \\
x_4 & x_2
                                  \end{array}\right),
\end{equation}
where $e_1=1,  e_2=-A_1 \overline{\tau_1},  e_3= A_2 \tau_2, e_4=e_2 e_3$.
The preimage of the lattice $L_N \subset \Vo$ is the $\Oc_{\k, 0}$-module
\begin{equation}
  \label{eq:Nf}
 \Nf := \gcd(\Oc_{\k, 0}\overline{\nf_1}, \Oc_{\k, 0}\nf_2) 
\end{equation}
with $\nf_i := N\Zb + NA_i\Zb \tau_i$ an $\Oc_i$-integral ideal of norm $N$.
Since $N$ is square-free and satisfies \eqref{eq:Heegner}, $\Nf$
is the unique $\Oc_{\k, 0}$-ideal such that
\begin{equation}
  \label{eq:nF}
\Nf \cap \Oc_{1} = \overline{\nf_1}\text{ and }\Nf \cap \Oc_{2} = \nf_2.
\end{equation}
Furthermore, the map \eqref{eq:WVo} identifies $\Sc(W(\Ab_F))$ and $\Sc(\Vo(\Ab))$, and $\Nf^\vee/\Nf \cong L_N^\vee/L_N$, 
where $\Nf^\vee$ is the $\Zb$-dual of $\Nf$ with respect to $Q_F$ \cite[section 2]{BL23}. 

In such a case, the associated torus $T$ over $\mathbb{Q}$ is  maximal in $\Ho$ and given by
$$T(R)=\left\{\left(t_1, t_2\right) \in\left(\k_1 \otimes_{\mathbb{Q}} R\right)^{\times} \times\left(\k_2 \otimes_{\mathbb{Q}} R\right)^{\times} \mid t_1 \bar{t}_1=t_2 \bar{t}_2\right\}$$
for any $\mathbb{Q}$-algebra $R$. The maps from $T$ to $\k$ and $\Ho$ are given by $\left(t_1, t_2\right) \rightarrow t_1 / {t}_2$ and $\iota(t_1, t_2) = (\iota_1(t_1), \iota_2(t_2))$, i.e.\
\begin{equation}
  \label{eq:TH-act}
  \iota(t) \cdot \phi(\alpha)
  = \phi((t_1/t_2) \alpha),~ t = (t_1, t_2) \in T, \alpha \in \k.
\end{equation}
This gives rise to the CM cycle
\begin{equation}
  \label{eq:ZW1}
  Z(W, z_0^\pm) = T(\mathbb{Q}) \backslash \{z^\pm_0\} \times T(\hat\Qb) / \iota^{-1}(K_0(N)) ,
\end{equation}
with $z_0^+ = (\tau_1, \tau_2) \in \Hb^2, z_0^- = \overline{z_0^+} \in (\Hb^-)^2$ \cite[section 2.2]{YY19}.
In addition, we have the CM cycle $Z(W, z_0^{\pm, \prime})$ with  $z_0^{+, \prime} = (-\overline{\tau_1}, \tau_2) \in \Hb^2, z_0^{-, \prime} = \overline{z_0^{+, \prime}} \in (\Hb^-)^2$.
Together, they give us the CM cycle on $\Xo_N$
\begin{equation}
  \label{eq:ZW}
  Z(W) =   Z(W, z_0^\pm) + Z(W, z_0^{\pm, \prime}),
\end{equation}
which is defined over $\Qb$.

Since $\iota_i$ are optimal for $i = 1, 2$, one has $\iota_i^{-1}\left(K_0(N)\right)=\hat{\mathcal{O}}_{i}^{\times}$.
Therefore, $\k_i^{\times} \backslash \hat\k_i^{\times} / \iota_i^{-1}\left(K_0(N)\right)$ is isomorphic to $\mathrm{Cl}\left(\k_{i}\right)$, the class group of $\k_i$. 
This gives rise to an injection (see Lemma 3.5 in \cite{YY19})
\begin{equation}
  \label{eq:p'}
  p': T(\mathbb{Q}) \backslash T(\hat\Qb) / \iota^{-1}(K_0(N)) \hookrightarrow \mathrm{Cl}\left(\k_{1}\right) \times \mathrm{Cl}\left(\k_{2}\right),~
\left[\left(t_1, t_2\right)\right] \to \left(\left[t_1\right],\left[t_2\right]\right). 
\end{equation}
By Proposition 3.3 in \cite{Li21}, it has the same image as
the natural map
\begin{equation}
  \label{eq:p''}
  p'': \operatorname{Gal}(H / \k) \hookrightarrow \operatorname{Gal}\left(H_{1} / \k_1\right) \times \operatorname{Gal}\left(H_{2} / \k_2\right)
 \cong  \mathrm{Cl}\left(\k_{1}\right) \times \mathrm{Cl}\left(\k_{2}\right) ,
\end{equation}
where $H := H_{1}H_{2}$ with $H_{i}$ the Hilbert class field of $\k_i$.
This gives us the parametrization (see \cite[Prop.\ 3.7]{YY19})
\begin{equation}
  \label{eq:ZW2}
  Z(W) = \sum_{\sigma \in \Gal(H/K)}
  (\tau_1^\sigma, \tau_2^\sigma)
  +     ((-\overline{\tau_1})^\sigma, (-\overline{\tau_2})^\sigma)
  +   ((-\overline{\tau_1})^\sigma, \tau_2^\sigma)
  +   (\tau_1^\sigma, (-\overline{\tau_2})^\sigma)
\end{equation}
By Lemma 3.2 in \cite{Li21}, the field $H_0 := H_1 \cap H_2$ is Galois over $\Qb$ with $\Gal(H_0/\Qb)$ an elementary 2-group.
In particular, it is the intersection of the genus subfields of $H_i$. 
When $D_0 := \gcd(D_1, D_2)$ is odd, it is given by
\begin{equation}
  \label{eq:H0}
  H_0 = H_0' :=  \Qb\lp\sqrt{(-1)^{(p-1)/2} p}: p \mid D_0\rp,
\end{equation}
which is disjoint from $F$ and $\k_i$.
The inclusion $H_0 \subset H_0'$ follows from $H_0$ being an elementary 2-group ramified only at primes dividing $D_0$, and $H_0 \supset H_0'$ follows from $H_0' \subset H_j$ for $j = 1, 2$  by genus theory.
So there exists $c_1, c_2 \in \Gal(H/\Qb)$ 
such that $c_i \mid_{\k_j}$ is complex conjugation when $i = j$ and identity otherwise.
%
Therefore, we can write
\begin{equation}
  \label{eq:CM-val}
  \Phi(Z(W)) = \sum_{\sigma \in \Gal(H/\Qb)} \Phi(\tau^\sigma_1, \tau^\sigma_2). 
\end{equation}
 for any function on $\Phi$ on $\Xo_N \cong Y_0(N)^2$.

When $\Phi = \Phi^{k-1}_{\mathrm{vv}(f)} = - G_{k, f}$ 
is the regularized theta lift associated with the harmonic Maass form $f$ as defined in section 2.4 of \cite{BLY22}, we have the following result.

\begin{theorem}
  \label{thm:CMformula}
In the notation above, for $f \in$ $H_{2-2k }(N)$ with $2 \nmid k$ we have
\begin{equation}
  \label{eq:CMval}
  -\frac{2}{|Z(W)|}
G_{k, f}(Z(W))
=
\mathrm{CT}\left(\tr^N_1\lp f^{+} \cdot\mathcal{C}_{k-1}\left(\Ec_N^{+}\rp\rp\rp
    - L'(0, \xi(f); D_1, D_2),
\end{equation}
where $\Ec_N^+$ is the holomorphic part of
the first derivative at $s = 0$ of the incoherent Hilbert Eisenstein series
$E_{\Nf, (1, 1)}(z, s)$ of parallel weight 1 defined in \eqref{eq:EP} associated with the $\infty_2$-neighbor of $W$ and $\Oc_{\k, 0}$-ideal $\Nf$, and
\begin{equation}
  \label{eq:LsG}
  L(s, \Go; D_1, D_2)
  :=
  \langle \Go(\cdot), \Cc_{k-1}(E_{\Nf, (1, 1)}(\cdot, \overline s)) \rangle_{\mathrm{Pet}}
\end{equation}
for any $\Go \in S_{2k}(N)$. 
  \end{theorem}

  \begin{remark}
    \label{rmk:N-dep}
    Though the Eisenstein series $E_{\Nf, (1, 1)}$ depends on $\Nf$, its diagonal restriction, or more generally image under $\Cc_{k-1}$, only depends on $N$ and $D_i$.
    Also, $E_{\Nf, (1, 1)}$ is the trivial component of the vector-valued incoherent Hilbert Eisenstein series
    $$
    E_{\Nf}(z, s) := (y_1y_2)^{-1/2}\sum_{\mu \in \Nf^\vee/\Nf}E(h_z, \cha(\hat \Nf + \mu)) \ef_\mu$$
    valued in $\Cb[\Nf^\vee/\Nf] \cong \Cb[L_N^\vee/L_N]$.
  \end{remark}

  \begin{proof}
    We apply Theorem 5.10 in \cite{BEY21} by taking $j = k-1$ even, in which case $Z^{k-1}(W)$ loc.\ cit.\ is $Z(W)$. 
     One can use
    \begin{equation}
      \label{eq:CTvv}
      \tr^N_1(f \cdot \Cc_{k-1} E_{\Nf, (1, 1)}))
      = 
      \langle \vv(f), \Cc_{k-1}E_{\Nf}\rangle
    \end{equation}
    to reduce the right hand side loc.\ cit.\ to that of \eqref{eq:CMval}. 
  \end{proof}

\section{Matching Sections}
In this section, we will deduce the following main result concerning the image of the Doi-Naganuma lift.
Throughout, $\chi_i, \chi, \k, \k_i$ are the same as in section \ref{subsec:Hecke} and $V, \Fc$ are the same as in section \ref{subsec:V}.
\begin{theorem}
  \label{thm:match}
  In the notations above, for any Siegel section $\Phi \in I(s, \chi)$ and function $C(s)$ real-analytic at $s = 0$ and $n \ge 1$,
  there exists
  a standard section $\Xi \in I(s, \chi_2)$  and $\varphi \in \Sc(V(\Ab))$ such that
  \begin{equation}
    \label{eq:match-global}
    \Fc(h, s; \Xi, \varphi) =
C(s)
\Phi(h, s) + O(s^n).
  \end{equation}
\end{theorem}
\begin{proof}[Proof of Theorem \ref{thm:main}]
  This follows immediately from Lemma \ref{lemma:unfold} and Theorem \ref{thm:match} above.
\end{proof}

\subsection{Archimedean Part}
\label{subsec:arch}
Let $\Phi_\infty = \Phi_\infty^{(k_1, k_2)} \in I(s, \chi_\infty)$ be the standard section, and right $K_\infty = \SO_2(\Rb)^2$-equivariant of weight $(k_1, k_2) \in \Zb^2$. Since $\chi_{\infty_i} = \sgn$ for both $i = 1, 2$, we have $k_i \equiv 1 \bmod{2}$.
Therefore, the quantities
\begin{equation}
  \label{eq:kl}
  k:= \frac{k_1+k_2}2,~
  l:= \frac{k_1-k_2}2
\end{equation}
are integers.
For our purpose, we will be mostly interested in the case $|k_1| = |k_2|$, i.e.\ either $k$ or $l$ is 0.
On the other hand, the result below holds for any odd $k_i$.

In the coordinate
$
V(\Rb) =
\{
\frac1{\sqrt{\aa}}
\smat{a}{\nu_1}{\nu_2 }{b}
: a, \nu_i, b \in \Rb
\}  $, we define the Schwartz function $\varphi^{ l,k}_\infty \in \Sc(V(\Rb))$ by
  \begin{equation}
    \label{eq:varphikl}
    \begin{split}
      \varphi_\infty^{l,k}(x)
      &:=
        p_+^l(x)
        p_-^{k} (x) 
        e^{-\frac\pi2
(        |(x, Z_+)|^2 +         |(x, Z_-)|^2)
        },\\
      p_+^{l}(x)
      &:=
        \begin{cases}
          \lp x, Z_+ \rp^l
= ((a+b)- i (\nu_1-\nu_2))^{l}  & l \ge 0,\\
          \lp x, \overline{Z_+} \rp^{-l}
= ((a+b)+ i (\nu_1-\nu_2))^{-l}& l \le 0,
        \end{cases}\\
      p_-^{k}(x)
        &:=
        \begin{cases}
          \lp x, \overline{Z_-} \rp^k
= ((a-b)+ i (\nu_1+\nu_2))^{k}          & k \ge 0,\\
          \lp x, {Z_-} \rp^{-k}
= ((a-b)- i (\nu_1+\nu_2))^{-k}           & k \le 0
        \end{cases},\\
      Z_+ &=X_+ + iY_+  := \aa^{-1/2} \smat{1}{-i}{i}{1},~
      Z_- =  X_- + iY_- := \aa^{-1/2} \smat{-1}{i}{i}{1}.
    \end{split}
  \end{equation}
  In addition, let $\Xi_\infty = \Phi^{|k| - |l|}_\infty \in I(s, \sgn)$ be the standard section satisfying
  \eqref{eq:keq}.
  Then we have the following result.

  \begin{lemma}
    \label{lemma:arch-match}
    Let $\Phi_\infty, \varphi_\infty = \varphi^{ l,k}_\infty, \Xi_\infty $ be as above. Then we have
    \begin{equation}
      \label{eq:arch-match}
      4\Fc_\infty(h, s; \Xi_\infty, \varphi_\infty)
      = \sqrt{\aa}^{-s-1} \Gamma_\Rb(s + 1 + |k| + |l|) \Phi_\infty(h, s)
    \end{equation}
    for all $h \in G(F \otimes \Rb)$ and $s \in \Cb$.
  \end{lemma}
  \begin{proof}
    It is easy to check that $\Fc_\infty(h, s) \in I(s, \chi_\infty)$ and right $\SO_2(\Rb)^2$-equivariant of weight $(k_1, k_2)$. Therefore, it suffices to prove \eqref{eq:arch-match} for $h = 1$.
    Using right $K_\infty$-equivariance of $\Xi_\infty$ and $\varphi_\infty$, we can simplify the expression for $\Fc_\infty$ in \eqref{eq:Frp} to
\begin{align*}
4\Fc_\infty(1, s)
  &=
    2    \int_0^{\infty} r^{s+1}
    \varphi^{ l,k}_\infty(x_0r) \frac{dr}r
    = \frac2{\sqrt{\aa}^{1+s}}     \int_0^{\infty} r^{s+1 + |k| + |l|} e^{-\pi r^2}    \frac{dr}r
    = \frac{ \Gamma_\Rb(s + 1 + |k| + |l|)}{\sqrt{\aa}^{1+s}}.
\end{align*}
\end{proof}


\subsection{Non-archimedean Part I}
\label{subsec:non-arch-I}
We begin by matching the unramified places.
\begin{proposition}
  \label{prop:unram}
  Suppose $p$ is unramified in $\k_1$ and $v_p(\aa) = 0$.
Let $\Phi_p \in I(s, \chi_p)$ be the right $G(F_p)$-invariant, standard section with $\Phi_p(1) = 1$,
  $\Xi_p  \in I(s,\chi_{2,p})$ be the right $\K_p$-invariant standard section with $\Xi_p(1) = 1$,
  and $\varphi_p$ be the characteristic function of $M_2(\Oc_p) \cap V_p \subset V_p$.
  Then
  \begin{equation}
    \label{eq:non-arch-match-ur}
    L_p(s+1, \chi_2)    \Fr_p(h, s)  =
    L_p(s+1, \chi)    \Phi_p(h, s)
  \end{equation}
  for all $h \in G(F_p)$ and $s \in \Cb$ with $\Re(s) > -1$.
  Here $L_p(s, \rho)$ is the local factor of $L(s, \rho)$ at the prime $p$ for $\rho = \chi, \chi_i$.
\end{proposition}
\begin{proof}
  Using right $\K_p$-invariance, it suffices to check \eqref{eq:non-arch-match-ur} for $h = 1$, in which case we can simplify $\Fc_p$ in \eqref{eq:Frp} to
  \begin{align*}
    \Fc_p(1, s)
    &=
    \int_{\Zb_p \backslash \{0\}} |a|_p^{s + 1} \chi_{1, p}(a)    d^\times a
      = (1-\chi_{1}(p)p^{-s-1})^{-1}
      =  L_p(s+1, \chi_1).
  \end{align*}
The claim then follows from \eqref{eq:Lfac}.
\end{proof}

For arbitrary standard section $\Phi_p$, we do not expect to match $L_p(s+1, \chi_1) \Phi_p$ on the nose as seen from the result below.
\begin{lemma}
  \label{lemma:standard}
  If $\Fc_p(h, s) = P(p^{-s-1}\chi_{1,p}(p))\Phi(h, s)$ for a standard section $\Phi \in I(s, \chi_p)$ and a non-trivial, non-invertible Laurent polynomial $P(X) \in \Cb[X, X^{-1}]$, then $\Phi(h, s)$ is independent of $h$ when $h \in H(\Zb_p)$.
\end{lemma}
\begin{proof}
  Since $P$ is non-trivial and non-invertible, we can write
 $
P(X) = X^N \lp \sum_{m = 0}^M c_m X^m\rp
$
with $N \in \Zb, M \ge 1$ and $c_0c_M \neq 0$.
Denote $\tilde\varphi = \Cc(\varphi) \in \Sc(V(\Qb_p))$.
For $X = \chi_{1, p}(p)p^{-s-1}$, we have
\begin{align*}
  P(X)
    \Fc_p(h, s)
  &= X^N \sum_{m= 0}^M
    c_m \Lambda_s(\tilde\varphi)(h/p^m)
  =\sum_{k\in\mathbb{Z}} X^{k+N}
    \sum_{m = 0}^M    c_m
\Lambda^*_0(\tilde\varphi)(h/p^{m+k}),
\end{align*}
where  $\Lambda^*_0(\tilde \varphi)(h) :=
    \int_{\Zb_p^{\times}}
    \tilde\varphi(h^{-1}x_0a) \chi_{1, p}(a)  d^\times a.
    $
For fixed $h \in H(\Zb_p)$, the function above is independent of $s$ if and only if
\begin{equation}
  \label{eq:induct-pol}
\sum_{m = 0}^M c_m \Lambda^*_0(\tilde \varphi)(h/p^{m+k}) = 0
\end{equation}
for all  $k \neq 0$.
For all  sufficiently positive $k$ (depending on $\tilde\varphi$ and $h$), the integral defining $\Lambda^*_0(\tilde\varphi)(h/p^k)$ becomes $\tilde\varphi(0) \int_{\Zb_p^\times} \chi_{1, a}(a)d^\times a$. On the other hand for all  sufficiently negative
$k$, we have  $\Lambda^*_0(\tilde\varphi)(h/p^k) = 0$ since the region of integration is outside the support of the Schwartz function $\tilde\varphi$.
In both cases, we can make remove the dependence on $h$ when it is in the compact set $H(\Zb_p)$.
A simple induction using \eqref{eq:induct-pol} then shows that $\Lambda^*_0(\tilde\varphi)(h)$ is independent of $h \in H(\Zb_p)$. Therefore,
$$
\Lambda_s(\tilde\varphi)(h) = \sum_{k \in \Zb} X^k \Lambda^*_0(h/p^k)
$$
is independent of $h$ when $h \in H(\Zb_p)$.
\end{proof}

For $i = 1, 2$, we define $\Cc_i : \Sc(V(\Qb_p)) \to \Sc(V(\Qb_p))$ by
\begin{equation}
  \label{eq:Cc1}
  \begin{split}
      \Cc_1(\varphi)
    &:=   \int_{J_p \cap N_p^-} \Xi_p(n^-(b))(\omega_p(n^-(b))\varphi)(\cdot)db,~\\
      \Cc_2(\varphi)
    &:=   \int_{ N_p^-} \Xi_p(n^-(b)w)(\omega_p(n^-(b)w)\varphi)(\cdot)db,
  \end{split}
\end{equation}
where $J_p$ and $N_p$ are defined in \eqref{eq:Jp}, and simply the expression of $\Fc_p(h, s)$ as follows.

\begin{lemma}
  \label{lemma:Fc1}
  For any $\Xi_p \in I(s, \chi_{2,p})$ and $\varphi_p \in \Sc(V(\Qb_p))$, we have
  \begin{equation}\label{decomposition Iwahori}
\begin{split}
  \Fr_p(h, s; \Xi_p, \varphi_p)
  &= \Lambda_s(\Cc_1(\varphi) + \Cc_2(\varphi)),
\end{split}
\end{equation}
where $\Lambda_s: S(V(\Qb_p)) \to I(s, \chi_p)$ is defined in \eqref{eq:Lambda}.
\end{lemma}
\begin{remark}
  In general we have $\Cc \neq \Cc_1+ \Cc_2$. But the difference lies in the kernel of $\Lambda_s$.
\end{remark}

\begin{proof}
  To simplify notation, we omit subscript $p$ from $\varphi_p, \Xi_p, J_p, N_p, M_p, N^-_p$.
  Using the decomposition \eqref{eq:Kp}, we can rewrite
  \begin{align*}
      \Fr_p(h, s)
    &  =     
      \int_{NM(J\cap N^-) \sqcup NMN^-w}\Xi(k)
\Lambda_s(\omega_p(k) \varphi)(h)dk.
  \end{align*}
By \eqref{eq:r-act}, the integral over $NM(J\cap N^-)$ simplifies to
  \begin{align*}
&      \int_{NM(J\cap N^-)}\Xi(k)
\Lambda_s(\omega_p(k) \varphi)(h)  dk
  =
  \int_{ J\cap N^-} \Xi(n^-) \Lambda_s(\varphi_p(k)\varphi)(h) dn^-
  =  \Lambda_s(\Cc_1(\varphi))(h).
  \end{align*}
  The integral over $NMN^-w$ simplifies analogously.
\end{proof}

\begin{remark}
  \label{rmk:Xip}
  Since $B(\Qb_p)(J_p \cap N_p^-)$ and $B(\Qb_p)N_p^-$ are open, disjoint with union being $G(\Qb_p)$ by the Iwasawa decomposition, we know that there exists unique $\Xi_p \in I(0, \chi_{2, p})$ such that
  \begin{equation}
    \label{eq:Xip}
    \Xi_p(k) =
    \begin{cases}
      1, & k \in J_p \cap N^-_p,\\
      0, & k \in N_p^-w.
    \end{cases}
  \end{equation}
  For such $\Xi_p$, $\Lambda_s(\Cc_2(\varphi))$ in \eqref{decomposition Iwahori} vanishes identically and
  \begin{equation}
    \label{eq:Fp1}
    \Fr_p(h, s; \Xi_p, \varphi_p)  =
    \Lambda_s(\Cc_1(\varphi))(h)
    = \int_{J_p \cap N^-_p}
\Lambda_s(\omega_p(n^{-}, h) \varphi_p) dn^-.
  \end{equation}
\end{remark}

\subsection{Non-archimedean part II}
\label{subsec:non-arch-II}
Despite of Lemma \ref{lemma:standard}, we can try to match $\Phi \in I(0, \chi)$ with $\Fc_p(h, 0; \Xi_p, \varphi_p)$ for suitable $\Xi_p$ and $\varphi_p$. To carry this out, we will utilize  the structure of $I(0, \chi_v)$ and Whittaker coefficients of $\Phi \in I(0, \chi_v)$ defined by
\begin{equation}
  \label{eq:Whittaker}
  W_m(\Phi) := \int_{F_v} \Phi(wn(b))\psi_v(-mb)db
\end{equation}
for $m \in F_v^\times$.
It has the property that $W_m(\Phi) = 0$ whenever $\Phi \in R(W_\alpha)$ and $m \not\in \alpha\Nm(\k_v^\times)$.

For $m=(m_v)_{v|p}\in F_p\cong \prod_{v|p}F_v$, the $m$-th Fourier coefficient of $\Fr_p (h,0)\in I(0, \chi_p)$ is
\begin{align*}
  W_m(\Fr_p (h,s))&= \prod_{v|p} W_{m_v}(\Fr_v(h,s)).
\end{align*}
We can now use this coefficient to give a surjectivity criterion of maps into $I(0, \chi_p)$.
\begin{proposition}
  \label{prop:surjective}
  Let $M$ be a $\Cb[H(\Qb_p)]$-module
  Then $\Fc \in \mathrm{Hom}_{\Cb[H(\Qb_p)]}(M,  I(0,\chi_p))$ is surjective if for every $\alpha \in F_p^\times /\Nm(\k_p^\times)$, there exists $m \in \alpha \Nm(\k_p^\times)$ and $\phi \in M$ such that $W_m(\Fc(\phi)) \neq 0$.
\end{proposition}

\begin{proof}
  Depending on the splitting behavior of $p$ in $F$ and $\k$, the representation $I(0, \chi_p)$ decomposes into direct sum of irreducible representations in the following way.
  $$
  I(0, \chi_p) =
  \begin{cases}
\oplus_{\alpha \in F_p^\times/\Nm(\k_p^\times)} R(W_\alpha),&  \text{$p$ is non-split in $F$,}\\
\oplus_{(\alpha, \alpha') \in F_p^\times/\Nm(\k_p^\times)} R(W_\alpha) \otimes R(W_{\alpha'}),& \text{$p$ is split in $F$.}
  \end{cases}
  $$
  The summands above are irreducible and pairwise non-isomorphic by the main theorems in \cite{KR92}.
  If $W_m(\Fc(\phi)) \neq 0$ for $m \in \alpha\Nm(\k_p^\times)$, then the projection of $\Fc(\phi)$ to the irreducible component corresponding to $\alpha \in F_p^\times/\Nm(\k_p^\times)$ is non-trivial. We are now done by the simple algebra lemma in \ref{lemma:surj}.
\end{proof}

Here is the simple algebra lemma mentioned in the proof above.
\begin{lemma}
  \label{lemma:surj}
  Let $R$ be a ring and $N_1, \dots, N_r$ pairwise non-isomorphic, simple $R$-modules.
  Denote $N:= \oplus_{1 \le j \le r} N_i$
and
$\pi_i: N \to {N}_i$ the natural projection for $1 \le i \le r$.
If $M \subset N$ is a submodule such that $\pi_i(M)$ is non-trivial for all $1 \le i \le r$, then $M = N$.
\end{lemma}

  This follows from Goursat's lemma for modules \cite{Kublik10}. For completeness, we include a proof here.

\begin{proof}
  We carry out the proof by induction. The case $r = 1$ is trivial since $N_1$ is simple. Suppose $r \ge 2$.
  For $1 \le i \le r$, denote ${N}^c_i := \oplus_{1 \le j \le r,~ j \neq i} N_i$ and  $\pi_i^c: N \to {N}^c_i$ the natural projection.
  By the inductive hypothesis, $\pi_i^c(M) = {N}^c_i$ for all $i$.
  Since $N_i$ is simple, $\ker \pi^c_i \cap M$ is either $N_i$ or trivial.
  If it is $N_i$ for one $i$, then we are done. Otherwise, $N_i^c \cong M/(M \cap \ker \pi^c_i) = M$ for all $i$. In particular we have $N_1^c  \cong N_2^c$, which is a contradiction since $N_i$'s are simple and pairwise non-isomorphic.
\end{proof}

Now we are ready to prove the following main result.

\begin{proposition}\label{prop:surj}
  The map $\Fc: I(0,\chi_{2,p})\times \Sc(V_p)\rightarrow I(0,\chi_p)$ sending $(\Xi, \varphi)$ to $\Fc_p(\cdot, 0; \Xi, \varphi)$
  is surjective.
\end{proposition}
\begin{proof}
  By Proposition \ref{prop:surjective}, it suffices to find
  for each $\alpha \in F_p^\times/\Nm(K_p^\times)$ a pair   $\Xi, \varphi$ and   $m \in F_p^\times$ representing $\alpha$ such that $W_m(\Fc(\Xi, \varphi)) \neq 0$.

When $F_p$ is a field, let $\Oc_p, \df_p^{-1}$ be its ring of integers and inverse different. When  $F_p = \Qb_p^2$, let them be $\Zb_p^2$.
  Then up to multiplicative constant, the Fourier transform of
$$\varphi =
\cha
\lp
\left(
\begin{array}{cc}
p^{-n}\Zb_p & \mathfrak{d}_p^{-1}\\
\mathfrak{d}_p^{-1} & 1+p^{n}\Zb_p
\end{array}
\right)\cap V(\Qb_p).
\rp \in \Sc(V(\Qb_p))$$
is
$$
 \psi_p(a) \cha
\lp \pmat{p^{-n}\Zb_p}{\Oc_p}{\Oc_p}{p^n\Zb_p}\cap V(\Qb_p)\rp(x)
$$
for $n \ge 1$ and $x= \smat{a}{\nu}{\nu'}{b} \in V(\Qb_p)$, which is invariant under $\omega(n)$ for $n \in w^{-1}J_pw \cap N_p$.
Therefore, $\omega(k)\varphi = \varphi$ for all $k \in J_p \cap N^-_p = w(w^{-1} J_pw \cap N_p)w^{-1}$, and $\Cc_1(\varphi) = \varphi$.

Let $\Xi_p$ be the same as in Remark \ref{rmk:Xip}, which then gives us
$$
\Fc(\Xi_p, \varphi) = \Lambda_0(\varphi)(h) = \int_{\Qb_p^\times} |a|_p \chi_{1, p}(a) \varphi(h^{-1}x_0a) d^\times a.
$$
Its $m$-th Whittaker coefficient is given by
\begin{align*}
  &W_m(\Fr(\Xi_p, \varphi))=
    \int_{F_p}  \int_{\mathbb{Q}_p^{\times}} |a|^{} \chi_{1,p}(a) \varphi
    \lp a
\left(\begin{array}{cc}
\Nm(b) & -b\\
-b^\prime & 1
\end{array}\right)
\rp d^{\times}a \psi_p(-\tr(mb)) db\\
  &=
    \int_{\df_p^{-1}} \int_{1+p^n\Zb_p} \chi_{1,p}(a)
    d^{\times}a \psi_p(-\tr(mb)) db
= \vol(1+p^n\Zb_p)    \int_{\df_p^{-1}} \psi_p(-\tr(mb)) db
\end{align*}
when $1+ p^n\Zb_p \subset \ker \chi_{1, p}$.
This is non-zero for all $m \in \Oc_p$. Since every $\alpha \in F_p^\times/\Nm(K_p^\times)$ has a representative in $\Oc_p$, we are done with the proof.
\end{proof}

\begin{remark}
  \label{rmk:sn}
  For any $n \ge 1$ and $\varphi \in \Sc(V(\Qb_p))$, \eqref{eq:r-act} gives us
  $$
  \Fc_p\lp h, s; \Xi_p, \lp 1 - \frac{\omega_p(m(p))}{\chi_{2, p}(p) p}\rp^n \varphi\rp
  = (1-p^{-s})^n \Fc_p(h, s; \Xi_p, \varphi)
  =  (s\log p )^n \Fc(\Xi_p, \varphi) + O(s^{n+1}).
  $$
  Also any linear combinations of Taylor coefficients of $\Phi \in I(s, \chi_p)$ is in $I(0, \chi_p)$.
  Therefore, for any function $C_p(s)$ real-analytic at $s = 0$ and $n \ge 1$, Proposition \ref{prop:surj} implies that we can find $\Xi_p, \varphi$ such that
  \begin{equation}
    \label{eq:match-sn}
    \Fc_p(h, s; \Xi_p, \varphi) = C_p(s) \Phi(h, s) + O(s^n).
  \end{equation}
\end{remark}

\subsection{Proof of Theorem \ref{thm:match}}
We are now ready to prove the global matching result in Theorem \ref{thm:match}.
Without loss of generality that $\Phi = \otimes_v \Phi_v$ and $\Phi_{\infty_i}$ has weight $k_i$.

Let $k, l$ and $\Xi_\infty = \Phi^{|k|-|l|}_\infty$ be as in section \ref{subsec:arch}.
At the places $p$ that are unramified in $\k$ and $\Phi_v$ is spherical for all $v \mid p$, we choose $\Xi_p$ spherical with $\Xi_p(1) = 1$ and $\varphi_p$ characteristic function of unimodular lattice in $V(\Qb_p)$. This include all but finitely many finite places, and we denote the set of the other the places by $S$, which is non-empty.
Fixed a $p_0 \in S$, apply Proposition \ref{prop:surj} to choose $\Xi_p$ and $\varphi_p$ for $p \in S$ such that
$$
\Fc_p(h, s; \Xi_p, \varphi_p) =
\begin{cases}
  \frac{4C(s)}{L^S(s+1, \chi_1)\Gamma_\Rb(s + 1 + |k|+|l|)} \Phi_p(h, s) + O(s^n),& p = p_0,\\
  \Phi_p(h, s) + O(s^n),& p \in S -\{ p_0\}.
\end{cases}
$$
Here $L^S(s, \chi_1) = \prod_{p \not\in S} L_p(s, \chi_1)$ is the partial $L$-function for $\chi_1$.
Using $\Fc = \prod_{p \le \infty} \Fc_p$ and Lemma \ref{lemma:arch-match}
gives us \eqref{eq:match-global}.

\section{Explicit Examples}
\label{sec:exs}

In this section, we will give an example demonstrating Theorems \ref{thm:match} and \ref{thm:main}.
When $\Phi^\infty$ is unramified, the matching in Theorem \ref{thm:main} is explicit and on the nose.
In general when $\Phi^\infty$ is slightly ramified, the explicit matching is exact everywhere, but is on the nose along the diagonal. 
First, we define and study certain invariant vectors in the Weil representation, which is of independent interest.
\subsection{Invariant vector}
\label{subsec:inv-vec}
To simplify the calculation of the map $\Cc$ in \eqref{eq:Cc}, we recall some invariant vectors of the Weil representation $\omega=\omega_{p}$.
We will use the Weil representation $\rho_A$ for a finite quadratic module $(A, Q)$.
This can be translated to $\omega$ via the pairing in \eqref{eq:C-pair}.
\subsubsection{Induction and Restriction Maps}
\label{subsub1}
Let $H \subset A$ be an isotropic subgroup, and $H^\perp \subset A$ its orthogonal complement, which in particular contains $H$.
The quotient $H^\perp/H$ is naturally a finite quadratic module of size $|A|/|H|^2$, and there is a linear induction map
$$
\uparrow^A_H: \Cb[H^\perp/H] \to \Cb[A],~ \ef_\mu \mapsto \sum_{\lambda \in H + \mu} \ef_{\lambda}
$$
that intertwines the Weil representations.
Using the bilinear form, we have its adjoint
$$
\downarrow^H_A: \Cb[A] \to \Cb[H^\perp/H],~ \ef_\lambda \mapsto \sum_{\mu \in H^\perp} \langle \ef_{\lambda}, \ef_\mu\rangle \ef_{\mu + H}.
$$

Using this map, one can produce invariant vectors inductively.
One  example is $A = A_0^+ \oplus A_0^-$, where the finite quadratic module $A_0^\pm = A_0$ has quadratic form $\pm Q_0$. Then $H_\pm = \{(\mu, \pm\mu) \in A: \mu \in A_0\} \subset A$ is totally isotropic, self-dual, and
\begin{equation}
  \label{eq:wA0}
\wf_\pm(A_0) := \uparrow^A_{H_\pm} \ef_0 = \sum_{\mu \in A_0} \ef_{(\mu, \pm\mu)} \in \Cb[A]
\end{equation}
are invariant vectors. If $A_0 = \Zb/p\Zb$ and $Q_0 =  \frac{ax^2}p$ with $p$ an odd prime and $a \in (\Zb/p\Zb)^\times$, then $\wf_\pm(A_0)$ spans $\Cb[A]^{\K_p}$.


\subsubsection{Fundamental Invariants}
\label{sec:rk3}
  \label{sec:u3}
%

  In \cite{MS22}, it was shown that any invariant vector are linear combinations of inductions of 6 types of fundamental invariants.
  One of them is from the trivial $A$.
  Here, we describe 3 of 5 non-trivial fundamental invariants using the orthogonal group.

  Let $K/\Qb_p$ be a ramified quadratic extension with valuation ring  $\Oc$, uniformizer $\varpi$, and different ideal $\df$.
  Denote
  \begin{equation}
    \label{eq:df-nm}
    \dfn = \Nm(\df)
    =
    \begin{cases}
      p & \text{ if } 2 \nmid p,\\
      4 \text{ or }8 & \text{ if } p =2.
    \end{cases}
  \end{equation}
Let $\chi$ be the quadratic character corresponding to $K/\Qb_p$, viewed as a character of $\Qb_p^\times$. Then $1 + \dfn \Zb_p \subset \ker \chi$.

  Consider the lattices
  \begin{equation}
    \label{eq:L0}
L_0 = \left\{
\dfn \pmat{a}{\nu}{\nu'}{b}: a, b \in \Zb_p, \nu \in \df^{-1}
  \right\}
\subset V_0 := L_0 \otimes \Qb_p
  \end{equation}
with quadratic form $Q_0 = \beta \cdot \det$, where $\dfn\beta \in \Zb_p^\times$.
The dual lattice $L_0^\vee$ is given by $M_2(\Oc) \cap V_0$.
Define the group $\Hc_{/\Zb_p}$ by
\begin{equation}
  \label{eq:Hc}
\Hc(R) := \left\{
h \in \GL_2(R \otimes \Oc): \det h \in R^\times
\right\}
\end{equation}
for any $\Zb_p$-algebra $R$.
Then $\Hc(\Zb_p)$ acts on $L_0$ by sending $\lambda \in L_0$ to $h \lambda  ({}^th')/\det(h)$,
which induces a map $\Hc(\Zb_p) \to \GSpin(L_0^\vee/L_0)$ that factors through $\Hc(\Zb_p/\dfn\Zb_p)$.
We slightly abuse notation and write $\Hc$ for $\Hc(\Zb_p)$.
View $\chi$ be a character of $\Qb_p^\times$ and hence of
$\Hc$ by composing with the determinant map.
For any $\mu \in L_0^\vee/L_0$, denote $\Hc_\mu \subset \Hc$ the subgroup fixing $L_0 + \mu$.
The following result is crucial in constructing invariant vectors in $\Sc_{L_0}$.
\begin{lemma}
  \label{lemma:key2}
  For any non-trivial $\mu  \in L_0^\vee/L_0 $, the followings are equivalent
  \begin{enumerate}
  \item $\mu$ is isotrpic.
  \item $L_0 + \mu \sim   \ell + L_0$ with $\ell  := \smat{1}{0}{0}{0} \in L_0$.
  \item  $\Hc_\mu \subset \ker \chi $.
  \end{enumerate}
\end{lemma}

\begin{proof}

  $(1) \Leftrightarrow (2)$:  if $\mu \in L_0^\vee/L_0 $ is isotropic, then there exists isotropic $\lambda \in L_0^\vee$ such that $\lambda \in L_0 + \mu$ by Hensel's lemma. In addition, $\lambda$ is $\Hc$-equivalent to $\ell$, which implies $(2)$ after modulo $L_0$.
  The converse is clear.

  $(2) \Rightarrow (3)$:
  The stabilizer of $L_0 + \ell$ is given by $          \Mc\Nc$, where
  \begin{equation}
    \label{eq:ell-stab}
    \begin{split}
    \Mc &:= \left\{ h \in \Hc: h \equiv\smat{ \alpha}{}{}{\alpha'} \bmod\df \text{ for } \alpha \in \Oc^\times\right\},~\\
    \Nc &:= \left\{ h \in \Hc: h \equiv \smat1{\beta}{}1\bmod\df \text{ for } \beta \in \Oc\right\}.
    \end{split}
  \end{equation}
  So it is clearly contained in $\ker \chi$. Since $\ker \chi  \subset \Hc$ is normal, we have $\Hc_\mu \subset \ker \chi$  for any $\mu$ equivalent to $L_0 + \ell$.

  $\neg (1) \Rightarrow \neg (3)$:
  Suppose first $2 \nmid p$. 
Given $\lambda = \smat{a}{\nu}{\nu'}b \in L_0^\vee - L_0$ such that its image in $L_0^\vee/L_0$ is not isotropic, i.e.\ $Q(\lambda) \not\in\Zb$.
If $a, b\in p\Zb_p$,
then
there exists $c \in \Zb_p^\times$ such that $\chi(c) = -1$ and $ac \equiv a \bmod \dfn, b/c \equiv b \bmod \dfn$, i.e. $\smat{c}{}{}1 \in \Hc_\lambda$ and $\Hc_\lambda \not\subset \ker \chi$.
Suppose $b \in \Zb_p^\times$.
Then we can replace $\lambda$ by $\smat{b}{}{}{1} \smat1{\beta}{}1 \cdot \lambda$ for suitable $\beta$ to suppose that $\nu = 0$ and $b = 1$. As $\lambda + L_0$ is not isotropic, we have $a \in \Zb_p^\times$. 
If $\chi(-a) = 1$, then we can find $h \in \Hc$ such that $h \cdot \lambda = \smat{0}{\tilde\nu}{\tilde\nu'}{0}$, and the argument before shows that $H_{h \cdot \lambda} \not\subset \ker \chi$.
If $\chi(a) = -1$, then $\smat{a}{}{}1 \smat{}{-1}{1}{}$ is in $H_\lambda$ but not $\ker \chi$.

Suppose now $\chi(-a) = -1$ and $\chi(a) = 1$. 
Then $\chi(-1) = -1$ and $\lambda$ defines a quadratic form on $\Oc^2$, making it isomorphic to the valuation ring $\Oc_M = \Oc[\sqrt{-a}]$ for the quadratic extension $M = K(\sqrt{-a})$ of $K$ with norm as the quadratic form.
Then the action of $\Oc_M^\times$ on $\Oc_M \cong \Oc^2$ gives an injection $\Oc_M^\times \hookrightarrow \GL_2(\Oc)$, and the subgroup  $U \subset \Oc_M^\times$ consisting of elements, whose norm to $K$ squares to $1$, is mapped to $H_\lambda$. 
Since $p > 2$, $K$ is the unique ramified extension of $\Qb_p$ and the norm map from $\Oc_M^\times$ to $\Oc^\times$ is surjective. In particular, there is $h \in \Oc_M^\times$ with norm $-1$. Its image in $\GL_2(\Oc)$ is in $H_\lambda$, but not in $\ker \chi$. 
%
This proves   $\neg (1) \Rightarrow \neg (3)$ for odd $p$, and the argument can be suitably adapted for $p = 2$. \footnote{Alternatively, one can check all cases by hand.}
\end{proof}


Now we define the following vector
\begin{equation}
    \label{eq:uK}
    \uf_K = \sum_{h \in \Hc/\Hc_{\ell}} \chi(h)^{-1} \cha(L_0 + h\cdot \ell) \in \Sc_{L_0}.
    \end{equation}
    As a consequence of Lemma \ref{lemma:key2}, we have the following result for rank 3 invariant vector. 
    \begin{proposition}
      \label{propLinv-3}
The $\chi$-isotypic subspace of $\Sc_{L_0}$, denoted by $\Sc_{L_0}^\chi$, is 1-dimensional and generated by the $\K_p$-invariant $\uf_K$.
    \end{proposition}

    \begin{proof}
      The subspace  $\Sc_{L_0}^\chi$ is generated by
      $$
      \sum_{h\in\Hc/(\Hc_\mu \cap \ker\chi)} \chi(h)^{-1} \cha(L_0 + h\cdot \mu)
      $$
    over $\mu \in L_0^\vee/L_0$, which vanishes if and only if $\Hc_\mu \not\subset \ker \chi$.
    By Lemma \ref{lemma:key2}, this happens precisely when $L_0 + \mu$ is not equivalent to $L_0 + \ell$. This shows that $\Sc_{L_0}^\chi$ is 1-dimensional and generated by $\uf_k$.
    Since the action of $\K_p$ and $\Hc$ commutes, $\K_p$ acts on $\Sc_{L_0}^\chi$ via a character $\rho$ that is trivial on $n(1) = \smat{1}{1}{}{1} \in \K_p$. Therefore, we have
    $$
\rho(w^4) = 1 = \rho((wn(1))^3) = \rho(w)^3,
$$
for $w = \smat{}{1}{-1}{} \in \K_p$,
which means $\rho(s) = 1$. So $\rho$ is trivial, and $\uf_K$ is $\K_p$-invariant.
    \end{proof}

    \begin{remark}
      It is easy to check that
    \begin{equation}
      \label{eq:u3-val}
      \uf_K(a\ell) = \uf_K \lp \smat{a}{}{}1 \cdot \ell\rp =
      \begin{cases}
        \chi(a) & a \in \Zb_p^\times,\\
        0 & \text{ otherwise}
      \end{cases}
    \end{equation}
    for all $a \in \Qb_p^\times$.
    The cases that $\dfn$ is an odd prime, 4, and 8 are the 3 types of non-trivial fundamental invariants in \cite{MS22}
    \end{remark}


    \begin{remark}
      \label{rmk:uK-odd}
      Suppose $p$ is odd.
Let $\pi: \GL_2(\Oc) \to \GL_2(\Fb_p)$ be the natural surjection and $\Hc_0 := \Hc \cap \ker(\pi)$. 
Then $\Hc_0 \subset \Hc_\ell$, $\Hc/\Hc_0 \cong \GL_2(\Fb_p)$ via $\pi$ and $\Hc_\ell/\Hc_0$ is congruent to $Z\cdot N \subset \GL_2(\Fb_p)$, with $Z$ the center and $N = \{n(b): b \in \Fb_p\}$ the unipotent.
The finite quadratic module $L_0^\vee/L_0$ is $\Syme(\Fb_p)$.
The set $\{h\cdot \ell: h \in \Hc/\Hc_\ell\}$ consists precisely the elements in $L_0^\vee/L_0$ with determinant 0, which are of the form $\epsilon \mu_{(a, c)} = h_{\epsilon, a, c} \cdot \ell$ with $\epsilon \in \Fb_p^\times$ and $\mu_{(a, c)} := \smat{a^2}{ac}{ac}{c^2}$ for $(a, c) \in \Fb_p^2 - \{(0, 0)\}$. Here
  $ h_{\epsilon, a, c}=\smat{a}{*}{c}{*}\smat{\epsilon}{}{}{1}$ with $\det\smat{a}{*}{c}{*}=1$. 
It is easy to check that $\chi(\det(h_{\epsilon, a, c})) = \chi(\epsilon)$.
So we can express
\begin{equation}
  \label{eq:uK2}
  \begin{split}
    \uf_K
    &  =
  \sum_{ \epsilon \in \Fb_p^\times} \chi(\epsilon)
  \lp
  \cha\lp L_0 + \epsilon \cdot \ell \rp
+  \sum_{j \in \Fb_p} 
  \cha\lp L_0 + \epsilon \cdot \smat{j^2}jj1\rp
      \rp\\
&= \sum_{\Lambda \in \Pb^1(\Fb_p), \epsilon \in \Fb_p^\times/(\Fb_p^\times)^2} \chi(\epsilon) \phi_{\Lambda, \epsilon}    ,
  \end{split}
\end{equation}
where $\phi_{\Lambda, \epsilon} := \sum_{v \in \Lambda} \cha(L_0 + \epsilon\cdot \mu_v)$
for $\Lambda \in \Pb^1(\Fb_p)$ and $\epsilon \in \Fb^\times_p$.
    \end{remark}
    
\begin{remark}
  \label{rmk:Hc}
  Let $K/\Qb$ be a quadratic extension with discriminant $\Delta$.
  Denote $\chi, \df, \Oc, \dfn, L_0$ the same as above, where $\Zb_p$ is replaced by $\Zb$, and $L_{0, p} := L_0 \otimes \Zb_p$ and $K_p := K \otimes \Qb_p$. 
  Then the vector
  \begin{equation}
    \label{eq:uK-global}
    \uf_K := \otimes_{p \nmid \dfn} \cha(L_{0, p}) \otimes_{p \mid \dfn} \uf_{K_p} \in \Sc_{L_0}^\chi
  \end{equation}
  is $\SL_2(\Zb_p)$-invariant.
  When $\dfn$ is odd, any  $\mu\in L_0^\vee/L_0$ has a representative $\smat{a}{c/2}{c/2}{b}   \in L_0^\vee$ with $a, b, c \in \Zb$ and
  \begin{equation}
    \label{eq:chid}
    \uf_K(\mu) =
    \begin{cases}
      \chi_\Delta([a, c, b])      & \text{if } \Delta \mid c^2 - 4ab \\
      0 & \text{otherwise,}
    \end{cases}
  \end{equation}
where $\chi_\Delta$ is the generalized genus character on binary quadratic forms defined in section I.2 in \cite{GKZ87}. 
\end{remark}

\subsection{Matching Example}
\label{subsec:ex-I}
Let $D_i < 0$ be distinct fundamental discriminants such that $D_1$ is odd. 
Denote $D_0 := \gcd(D_1, D_2) \ge 1, D := D_1D_2$ and $\k, F, \k_i$ as in the introduction with discriminants $D_\k, D_F, D_i$. 
We rescale the quadratic forms
\footnote{It's easy to modify the argument in Section 3 to get the same Theorem \ref{thm:match}.}
on $\k_2$ and $V$ to be $D_2 \cdot \Nm$ and $|D_1|\cdot \det$, i.e.\ $\aa = |D_1|$ in \eqref{eq:V}, and consider the quadratic space
\begin{equation}
  \label{eq:tV}
  \tV := E_2 \oplus V = E_2 \oplus \Qb^2  \oplus F,~ \tQ(\alpha, a, b, \nu) = D_2\Nm_{\k_2}(\alpha) - D_1(ab - \Nm_{F}(\nu)).
\end{equation}

Suppose $N \in \Nb$ is square-free and satisfies condition \eqref{eq:Heegner}. In particular, all prime factors of $N$ are split in $\k_i$.
Using the invariant vectors in section \ref{sec:rk3}, we will first construct an invariant vector in $\tv^\infty_N = \otimes_{p < \infty} \tv_{N,p} \in \Sc(\hat\k_2 \oplus \hat V)$ as follows.

For any $\nf_2 \subset \Oc_2$ with $\Nm(\nf_2) = N$, take the lattice
\begin{equation}
  \label{eq:latt-M}
  \tM := \nf_2 \df_{\k_2}^{-1} \oplus M \subset \tV,~
  M   :=
  \Zb \oplus N\Zb \oplus \df_{F, 1}^{-1}
\subset M^\vee =  N^{-1}   \Zb \oplus \Zb \oplus (D_1 \df_{F, 2})^{-1}
 \subset V,
\end{equation}
where $\df_{F, j} \subset \Oc_F$ is the unique ideal such that $D_{F, j} := \Nm(\df_{F, j}) = \gcd(D_j, D_F) $ for $j = 1, 2$.
Note that
\begin{equation}
  \label{eq:Ds}
  D_F = \frac{D}{D_0^2},~
  D_{F, j} = \frac{D_F D_0}{|D_{j'}|} = \sqrt{\frac{D_jD_F}{D_{j'}}}
  = \frac{|D_{j}|}{D_0}
      \end{equation}
 with $j' := 3-j$,   since $D_0$ is odd. 
The dual of $\df_{F, 1}^{-1}$ with respect to $|D_1|\cdot \Nm$ is $\frac{1}{|D_1|} \df_F^{-1}\df_{F, 1} = \frac{1}{|D_1|} \df_{F, 2}^{-1}$.
Using $2 \nmid D_1$, we can explicitly describe $\tM^\vee$ as
\begin{equation}
  \label{eq:tMv}
  \begin{split}
  \{ (\alpha, a, b, c, d):
&  \alpha \in \nf_2, a, b, c, d \in \Zb, cD_{F} \equiv d \bmod 2  \}
\cong      \tM^\vee  ,\\
(\alpha, a, b, c, d) &\mapsto
    \lp \frac{\alpha}{ND_2}, \frac{a}{ND_1},\frac{b}{D_1} , \frac{d D_{F, 1} + c \sqrt{D_F}}{2 D_1 \sqrt{D_F}}\rp
  \end{split}
\end{equation}
and will sometimes use the tuple $(\alpha, a, b, c, d)$ to represent elements in $\tV$.
The quadratic form is given by
\begin{equation}
  \label{eq:tMv-Q}
 Q(\alpha, a, b, c, d) = \frac{\Nm_{\k_2}(\alpha)}{N^2D_2} - \frac{ab}{ND_1} + \frac{c^2}{4D_1} - \frac{d^2}{4D_2}. 
\end{equation}

Let $H := \langle \mu_j + \tM: j = 1, 2, 3 \rangle \subset \tM^\vee/\tM$ be the isotropic subgroup generated by
\begin{equation}
  \label{eq:ell}
  \begin{split}
      \mu_1 &= \lp \frac{N}{D_2} , 0, 0,  \frac{N}{\sqrt{D}}  \rp = (N^2,0,0,0,2N),\\
    \mu_2 &= \lp \frac{\nu}N , \frac1N, N,  0  \rp = (\nu D_2,D_1,ND_1,0,0),\\
  \mu_3 &= \lp D_1 , 1, D_2 \tr(\nu),  0  \rp = (ND,ND_1,D\tr(\nu),0,0), 
  \end{split}
\end{equation}
with $\nu \in \nf_2^2$ such that $\nu \not\in N\Oc_2$. 
It is easy to check that $(\alpha, a, b, c, d) \in \tM^\vee$ lies in $H^\perp + \tM$ if and only if
\begin{equation}
  \label{eq:Hperp-cond}
  \tr(\alpha) - Nd \in D_2\Zb,~
  \frac{\tr(\alpha \overline \nu)}N - b \in N\Zb,~
  D_1\tr(\alpha ) - D_2 \tr(\nu) a  \in N\Zb.
\end{equation}
In particular if $\alpha \in \tM$, then $(\alpha, a, b, c, d) \in \tV$ is in $H^\perp + \tM$ if and only if
\begin{equation}
  \label{eq:Hperp-cond0}
  d \in D_2\Zb,~
  a, b \in N\Zb,~
    c \in \frac1{(2, D_2)}(2\Zb + d).
  \end{equation}
  The same holds after tensoring with $\hat\Zb$.
The following lemma furthermore describes $H^\perp \subset \tM^\vee/\tM$.

\begin{lemma}
  \label{lemma:Hperp}
  In the notation above with $D_1, D_2 < 0$ fundamental discriminant. 
  Suppose $D_1$ is odd.  
  Then the finite quadratic module $L_0^\vee/L_0
\cong (\Zb/|D_1|\Zb)^3  $ from \eqref{eq:L0}   and Remark \ref{rmk:Hc} by taking $K = \k_1$, $d = |D_1|$ and $\beta = 1/|D_1|$.
  embeds into $\tM^\vee/\tM$ via

  $$
  \iota
\lp \pmat{a'}{c'}{c'}{Nb'} \rp
:= \lp 0, \frac{a'}{D_1}, \frac{Nb'}{D_1}, \frac{c'}{D_1}\rp + \tM.
$$
Furthermore, $\mathrm{im}(\iota) \subset M^\vee/M \subset \tM^\vee/\tV$ and $H^\perp = H \oplus \mathrm{im}(\iota)$. 
\end{lemma}


\begin{proof}
    It is easy to check that $\iota$ is an injective isometry, and its image is orthogonal to $H$ while intersects $H$ trivially.
  Therefore $H \oplus \mathrm{im}(\iota) \subset H^\perp$. 
  Since $|\tM^\vee/\tM| = N^4|D_1^3D_2^2|$ and $|H| = N^2|D_2|$, we also have
  $$
  |H^\perp| = |\tM^\vee/\tM|/|H| = N^2 D^3_1D_2 = |H| \cdot |\mathrm{im}(\iota)|.
  $$
  This finishes the proof. 
\end{proof}

We now apply this lemma 
to define the $\SL_2(\hat\Zb)$-invariant vector
\begin{equation}
  \label{eq:tvf}
  \tv_N^\infty = \otimes_{p < \infty} \tv_{N, p} := \uparrow^{A}_H \iota(\uf_{K_1}) \in \Sc_{\tM} 
\subset  \Sc(\hat\k_2 \oplus \hat V),
\end{equation}
with $A := \tM^\vee/\tM$ and $\uf_{\k_1}$ defined in \eqref{eq:uK-global}. 
For $\alpha \in N\nf_2 \df_2^{} \otimes\hat\Zb, \mu' \in L_0^\vee/L_0, d \in \hat\Qb$, we have
  \begin{equation}
    \label{eq:tvp-val}
    \tv^\infty_{N}(\alpha, \iota(\mu'), d) =
    \uf_{\k_1}(\mu')
\cha_{D_2\hat\Zb}(d)   .
      \end{equation}


In section\ref{subsec:CM}, 
we saw that  fixed CM points $\tau_i \in \Hb$ with discriminant $D_i$ give identification $\Sc(W(\Ab_F)) \cong\Sc(\Vo(\Ab))$ for $F$-quadratic space $W = \k$. In $\Sc(\hat\Vo)$, we have the $K_0(N)$-invariant Schwartz function $\cha(\hat L_N)$. 
Its image in $\Sc(W(\hat F))$ is $\cha(\hat\Nf)$ for a certain $\Oc_{\k, 0}$-ideal $\Nf$.
Using the invariant vector $\tv_N^\infty$, we can match the Hilbert Eisenstein series constructed from $\cha(\hat \Nf)$. 

\begin{proposition}\label{prop:ex1}
  Let $D_i < 0$ be fundamental discriminants with $D_0 := \gcd(D_1, D_2)$ odd.
  For square-free $N\in \Nb$ satisfying \eqref{eq:Heegner} and odd $k \ge 1$, we have
            \begin{equation}
              \label{eq:ex-main}
              \begin{split}
  \sum_{N' \mid N}
\frac{\zeta_N(s+1) N'/N }{\zeta_{N'}(1)\zeta_{N/N'}(s)}
\Fc(h, s; \tv^k_{N'})
&              = 4^{-1} \Gamma_\Rb(s + 1 + k) L^{}(s+1, \chi_1) |D_1|^{\frac{s+1}{2}}
  \Phi_{N}^k(h, s)
              \end{split}
            \end{equation}
            for all $h \in B(\Ab_F)K(D_0N)K_\infty \subset
            H(\Ab)$. 
            On the left hand side, $\zeta_N(s) := \prod_{p \mid N}(1-p^{-s})^{-1}$ is the partial zeta function, $\Fc(h, s; \tv) \in I(s, \chi)$ is defined in \eqref{eq:sec} (see also Remark \ref{rmk:1-2}), $\tv^k_N =
\tv_N^\infty \otimes \tv_\infty^k
\in \Sc(\tV(\Ab))$ with $\tv_N^\infty$
defined in \eqref{eq:tvf}, $\tv_\infty^k = \varphi_\infty^k \otimes \varphi_\infty^{( 0,k)} \in \Sc(\k_2 \otimes \Rb) \otimes \Sc(V(\Rb))$ with $\lambda(\tv_\infty^k) = \Phi^k_\infty$ and $\varphi^{( 0,k)}_\infty$ defined in \eqref{eq:varphikl}.
On the right hand side, $\Phi_N^k = \Phi_N^\infty \otimes \Phi_\infty^{(k, k)}$ with $\Phi_N^\infty = \otimes_{v < \infty}\Phi_{N, v} = \lambda(\cha(\hat\Nf)) \in I(s, \chi)$, where  $\Nf$ is any integral $\Oc_{\k, 0}$-ideal with norm $N$.
Also, $K(D_0N; F) \subset \SL_2(\widehat \Oc_F)$ is the preimage of $\SL_2(\Zb/(D_0N)\Zb) \subset \SL_2(\Oc_F/(D_0N))$ under the natural surjection modulo $D_0N$.
In particular, $K(D_0N; F)$ contains the subgroup $\SL_2(\hat\Zb) \subset \SL_2(\hat\Oc_F)$.
\end{proposition}

\begin{remark}
  \label{rmk:match-exact}
  If $D_0 = N = 1$, then
$K(D_0N)= \SL_2(\widehat\Oc_F)$ and
\eqref{eq:ex-main} holds for all $h \in B(\Ab_F)
 \SL_2(\widehat\Oc_F) K_\infty =  \SL_2(\Ab_F) = H(\Ab)$. 
\end{remark}

\begin{proof}
For the infinite places, applying Lemma \ref{lemma:arch-match} with $d = |D_1|$ gives us
\[
  4\Fc_\infty(h, s; \tv^k_\infty) =
  4\Fc_\infty(h, s; \Phi^k_\infty, \varphi_\infty^{( 0,k)}) =\frac{1}{|D_1|^{\frac{s+1}{2}}} \Gamma_\Rb(s + 1 + k) \Phi^{(k, k)}_\infty(h, s). \]
At a finite place $p$, denote $N_p := \gcd(N, p)$.
Since $N$ is square-free, it suffices to show that
\begin{equation}
  \label{eq:match-ex}
  \sum_{N' \mid N_p}
\frac{\zeta_{N_p}(s+1)N'/N_p}{\zeta_{N'}(1)\zeta_{N_p/N'}(s)}
  \Fc_p(k_p, s; \tv_{N', p})
  = |D_1|_p^{-s-1}
    L_{p}(s+1, \chi_1)
\prod_{v \mid p} \Phi_{N, v}(k_v, s)
\end{equation}
for all $k_p = (k_v)_{v \mid p} \in K(D_0N; F)_p$.

When $p \nmid D_0N$, 
both  $\tv_{N, p}$ and the section $\prod_{v \mid p} \Phi_v = \lambda(\cha(\Nf_p))$ are right $\SL_2(\Oc_p)$-invariant,
so we only need to check \eqref{eq:match-ex} for $k_p = 1$, where the right hand side is then just $|D_1|_p^{-s-1}L_p(s+1, \chi_1)$.
Furthermore, $\tv_{N, p}$ is $\K_p$-invariant by its construction.
So the left hand side simplifies as 
$$
\Fc_p(1, s; \tv_{N,p})
=
\tilde\Lambda_s(\tv_{N,p})(1)
=
 \int_{\Qb_p^\times}|a|_p^{s+1} \chi_{1,p}(a)\tv_{N,p}(0, a, 0, 0)d^\times a
 $$
by \eqref{eq:tFcp}.
Applying equations \eqref{eq:tvf} and \eqref{eq:uK2} shows that this is precisely $|D_1|_p^{-s-1}L_p(s+1, \chi_1)$.

When $p \mid ND_0$, we only know that the image $\cha(\Nf_p)$ in $\Sc(\widehat \Vo)$ is $\cha(\widehat L_N)$, which is invariant with respect to the subgroup $K_0(N)$ of $\SL_2(\Zb_p) \subset \SL_2(\Oc_p)$.
Also, $\cha(\Nf_p)$ is invariant with respect to the kernel of $\SL_2(\Oc_p) \to \SL_2(\Oc_F/p)$.
For $p \mid D_0$, the same argument above by checking $k_p = 1$ proves \eqref{eq:match-ex} when $k_p \in K(D_0N; F)_p \cap \SL_2(\Oc_p)$.

For $p \mid N$, we also need to check $k_p = n(j)w \in \K_p \subset \SL_2(\Oc_p)$ for $1 \le j \le p$.
On the one hand, $(L_N^\vee/L_N) \cong (\Zb/N\Zb)^2$ is just a scaled hyperbolic plane, and
we have 
\begin{align*}
\Phi_{N, p}(k_p) =
  (\omega(k_p) &\cha(L_{N, p}))(0)
  =\langle \omega(n(j)w)\mathfrak{e}_0,\mathfrak{e}_0\rangle
  = \left\langle \frac{\omega(n(j))}N\sum_{\mu \in (\Zb/N\Zb)^2} \ef_\mu, \ef_0\right\rangle = \frac1N,
\end{align*}
where $\ef_\mu = \cha(\mu + (N\Zb)^2 )$.
On the other hand, the $\K_p$-invariance of $\tv_{N,p}$ again gives us
\begin{align*}
  \Fc_p(k_p, s; \tv_{N,p})
&=
 \int_{\Qb_p^\times}|a|_p^{s+1} \chi_{1,p}(a)\tv_{N, p}(0, k_p^{-1}(a, 0, 0))d^\times a\\
&=
  \int_{\Qb_p^\times}|a|_p^{s+1} \chi_{1,p}(a)\tv_{N, p}(0, a(j^2, 1, -j))d^\times a
  = p^{-s-1} \chi_{1, p}(p) L_p(s+1, \chi_1),
\end{align*}
and similarly $  \Fc_p(k_p, s; \tv_{1,p}) = L_p(s+1, \chi_1)$.
Since $N$ satisfies \eqref{eq:Heegner}, we have $\zeta_p(s) = L_p(s, \chi_1)$. 
Substituting these into \eqref{eq:match-ex} verifies it, and completes the proof.
\end{proof}

\begin{proof}[Proof of Theorem \ref{thm:ex}]
    Combining Proposition \ref{prop:ex1} with Lemma \ref{lemma:unfold} above, we obtain 
    \begin{equation}
      \label{eq:matchN}
    \sum_{N' \mid N}
\frac{(N'/N) \zeta_N(s+1)}{\zeta_{N'}(1)\zeta_{N/N'}(s)}
I(h, \tv^k_{N'}, s)
= 4^{-1} \Gamma_\Rb(s + 1 + k) L^{}(s+1, \chi_1) |D_1|^{\frac{s+1}{2}}
E(h,\Phi_{N}^k, s)
            \end{equation}
            for $h \in  B(\Ab_F)K(D_0N)K_\infty \subset H(\Ab)$. 
The case $D_0 = N = k = 1$ reduces to Theorem \ref{thm:ex}. 
\end{proof}

\begin{remark}
Note that the quadratic form on $E_2$ is $D_2 \Nm_{E_2/\Qb}$. $k_1=k_2=1$ implies $k=1,~ l=0$ in Lemma \ref{lemma:arch-match}, and hence the elliptic Eisenstein series associated to $\otimes_{p<\infty} \Xi_p \otimes \Xi_\infty^1$ is incoherent.
\end{remark}

\section{CM values in higher level}
\label{sec:last}
In this section, we will use the theta lift expression of incoherent Hilbert Eisenstein series to explicate the generalized Rankin-Selberg $L$-series appearing in Theorem 1.1 of \cite{BKY12}, and its higher weight analogue, for the CM cycle in section \ref{subsec:CM}.
This turns out to involve Fourier coefficients of half-integral weight modular forms.
Using bounds for such coefficients, along with equidistribution results for the CM cycle in section \ref{subsec:CM}, we will prove a higher level generalization of \cite{Li21}. 

\subsection{Twisted Shintani Lift of Eigenforms}

For $k \in \Nb$, let $\Delta$ be a fundamental discriminant such that $(-1)^k\Delta > 0$. 
Consider the quadratic spaces
\begin{equation}
  \label{eq:V0}
  (V_{0}, Q_0)=(\Syme(\mathbb{Q})^{},|\Delta|\cdot \det).
\end{equation}
The group $\PGL_2$ acts on $V_0$ via $v \mapsto h v h^t/\det(h)$, which identifies it with $H_0 := \GSpin(V_0)$.
Denote $\omega_0 = \omega_{V_0}$ the associated Weil representation of $G'_\Ab \times H_0(\Ab)$ on $\Sc(V_0(\Ab))$. 
For $N_1, N_2 \in \Nb$ such that $N_1, N_2, 2\Delta$ are pairwise  co-prime, define the Schwartz function $\varphi^k_{  \Delta, N_1, N_2} = \otimes_{p < \infty} \varphi_{  \Delta, N_1, N_2, p} \otimes \varphi_{0, \infty}^k \in \Sc(V_0(\Ab))$ by
\begin{equation}
  \label{eq:varphi0N}
  \begin{split}
    \varphi_{  \Delta, N_1, N_2, p}
    &:=
      \begin{cases}
        \cha\lp\Syme(\frac12 \Zb_p) \cap \smat{N_2^{-1}\Zb_p}{\frac12\Zb_p}{\frac12\Zb_p}{N_1N_2\Zb_p} \rp & p \nmid \Delta,\\
\uf_{\Qb_p(\sqrt\Delta)} & p \mid \Delta,        
      \end{cases}\\
    \varphi_{0, \infty}^k(x)
    &:= 
          \lp x, \overline{Z_-} \rp^k
      e^{-\frac\pi2(|(x, Z_-)|^2 + (x, X_+)^2)      } \in \Sc(V_0(\Rb)),
  \end{split}
\end{equation}
where $Z_\pm$ is defined in \eqref{eq:varphikl} with $\aa = |\Delta|$ and $\uf_K$ is the invariant vector defined in \eqref{eq:uK}.
We can also view $\varphi_{  \Delta, N_1, N_2}^\infty := \otimes_{p < \infty} \varphi_{  \Delta, N_1, N_2, p}$ as in $\Sc_{\Lo}$ for the lattice
\begin{equation}
  \label{eq:L0N}
  \Lo
= \Lo_{N_1, N_2}
  :=
  \left\{
    \pmat{a}ccb: a \in
    N_2^{-1}\Zb,~ c \in\Zb,~ b \in N_1N_2\Zb
\right\},
\end{equation}
with support in $(\tfrac1{2|\Delta|} \Lo) \cap \Lo^\vee$.

Suppose $N_1, N_2$ are co-prime and square-free. 
It is easy to verify that
\begin{equation}
  \label{eq:AL-op}
  \begin{split}
    \omega_0\lp \smat{}{-1}{N_1'}{}_{N_1'} \rp
  \varphi_{  \Delta, N_1, N_2}^\infty
  &= \varphi_{  \Delta, N_1, N_2}^\infty,\\
    \omega_0\lp w_{N_2'} \smat{N_2'}{}{}{1}_{N_2'} \rp
 \varphi_{  \Delta, N_1, N_2}^\infty
&=    \omega_0\lp w_{N_2'} \rp
 \varphi_{  \Delta, N_1, N_2/N_2'}^\infty
 = \varphi_{  \Delta, N_1, N_2/N_2'}^\infty
  \end{split}
\end{equation}
for any $N'_{j} \mid N_j$ satisfying $\gcd(N'_j, N_j/N'_j) = 1$. 
Here $w = \smat{}{-1}1{} \in H_0(\Qb)$ is an involution and the notation $w_M \in H_0(\hat\Qb)$ is defined in the beginning of section \ref{sec:prelim}. 


We can now define the twisted theta function
\begin{equation}
  \label{eq:theta-Delta}
  \begin{split}
      \theta^{(k)}_{\Delta, N_1, N_2}(\tau, z)
    &:=
      \sum_{\mu \in L^\vee /L}
\varphi^{\infty}_{  \Delta, N_1 N_2}      (\mu)
      \theta^{(k)}_{\Lo+\mu}(\tau, z),\\
    \theta^{(k)}_{\Lo+\mu}(\tau, z)
    &=
      \frac{v^{k+1}}{      y^{2k}}
      \sqrt{|\Delta|}^k
      \sum_{X = \smat{a}ccb \in \Lo + \mu}
      (b\bar z^2 - 2c \bar z + a)^k
      \ebf(Q_0(X) u) e^{-2\pi v Q_{0, z}(X)},
  \end{split}
\end{equation}
where $Q_{0, z}(X) :=
\frac1{4|\Delta| y^2}
\lp {|(X , \smat{|z|^2}xx1)|^2}
+ {|(X , \smat{z^2}zz1)|^2}\rp
$
is the majorant.
It is modular in $\tau$, resp.\ $z$, of weight $-k-1/2$, resp.\ $2k$. 
Using Remark \ref{rmk:Hc}, we see that the theta kernel $\theta^{(k)}_{\Lo + \mu}$ is $\sqrt{|\Delta|}^{1-k}v^{k-1/2}$ times the $\ef_\mu$-component of the theta kernel $\Theta_{Sh}(\tau, z)$ in section 4 of \cite{AS21} for the lattice $\Lo$.
  \footnote{Note our parameter $k$ is $k+1$ in \cite{AS21}.}

We can also write \cite[Eq.\ (1.37)]{Kudla03}
\begin{equation}
  \label{eq:theta-Delta1}
      \theta^{(k)}_{\Delta, N_1, N_2}(\tau, z)
=   \sqrt{v}^{k+1/2} y^{-k}
    \theta_{V_0}(g'_\tau, h_z, \varphi^k_{  \Delta, N_1, N_2}),
\end{equation}
where $g'_z = [g_z, 1] \in G'_\Rb$, resp.\ $h_z \in \mathrm{PSL}_2(\Rb) \subset H_0(\Rb)$, is the preimage, resp.\ image, of $g_z \in \SL_2(\Rb)$.
We now state the Shintani lift of weight $2k$ eigenform via $\theta_{\Delta, N_1, N_2}^{(k)}$. 

\begin{proposition}
  \label{prop:Shintani}
  For $k, \Delta$ as above, let $N_0, N_1, N_2, N_3 \in \Nb$ be integers such that
 $N := \mathrm{lcm}(N_0N_3, N_1N_2)$ is square-free and co-prime to $2\Delta$, 
and $N_0$ satisfies condition \eqref{eq:Heegner} with respect to $\Delta$. 
For a normalized newform $\Goo = \sum_{m \ge 1} a_{\Goo}(m)q^m \in S_{2k}(N_0)$, we have
  \begin{equation}
    \label{eq:Shintani-eigenform}
    \begin{split}
      \int_{\Gamma_0(N) \backslash \Hb}
      &\Goo(N_3 z)
      \overline{\theta^{(k)}_{\Delta, N_1, N_2}(\tau, z
      )} y^{2k}d\mu(z)
      =
v^{k-1/2}        |\Delta|^{(1-k)/2}
\overline{c_{\tGoo}(|\Delta|)}
        2^{-k} (-1)^{k-1+\lfloor k/2 \rfloor}
        \frac{|\Goo|_{\mathrm{Pet}}^2}{6|\tGoo|_{\mathrm{Pet}}^2}
      \\
      &       \times
N_3^{-k}\sigma_1(N/N')
        \lp\frac{N'}{N_0}\rp^{1-k}
\sum_{c \mid (N/N_0)} \sum_{\rr \mid c}
(-1)^{\omega(\rr)}        a_{\Goo}\lp \frac{N'}{N_0 c}\rp \lp U_{\rr^2}
        \tGoo\rp(\tau/4)
    \end{split}
    \end{equation}
    if $N_0 \mid N_1$ and $\Goo$ is invariant with respect to the Fricke involution $\Wc_{N_0}$. 
    Otherwise this integral is identically 0.
Here $N' := N/\gcd(N_2, N_3)$, $\sigma_1$ is the sum of divisor function, $\omega(r)$ is the number of prime divisors of $r$, and $\tGoo(\tau) = \sum_{m \ge 1} c_{\tGoo}(m)q^m \in S_{k+1/2}^{+, \mathrm{new}}(4N_0)$ is the newform associated with $\Goo$ under the Shimura correspondence. 
\end{proposition}

\begin{proof}
  Note that $N$ being square-free implies $N_j$ is square-free and $\gcd(N_j, N_{3-j}) = 1$ for all $0 \le j \le 3$. 
  To start, we define
  \begin{equation}
    \label{eq:IDel}
      \begin{split}
  I_{\Delta, N_1, N_2, N_3}(g_\tau', \Goo^\#)
&:=
\sqrt{v}^{-k-1/2}
      \int_{\Gamma_0(N) \backslash \Hb} \Goo(N_3 z)
      \overline{\theta^{(k)}_{\Delta, N_1, N_2}(\tau, z
      )} y^{2k}d\mu(z)\\
&=
\frac{\vol(\Gamma_0(N)\backslash \Hb)}{\vol([H_0])}
\int_{[H_0]} (V_{N_3}\Goo)^\#(h)
\overline{  \theta_{V_0}(g_\tau', h, \varphi^k_{  \Delta, N_1, N_2})}
dh.
\end{split}
  \end{equation}
    If there is $p \mid N_0$ such that $p \nmid N_1$,
    we can apply \eqref{eq:AL-op} with $N_2' = p' :=\gcd(N_2, p)$
    to obtain
    \begin{align*}
      I _{\Delta, N_1, N_2, N_3}(g', \Goo^\#)
      &=
        I _{\Delta, N_1, N'_2, N_3}\lp g',\Wc_{p'}\Goo^\#\rp
=         I _{\Delta, N_1, N'_2, N_3}\lp g',\tr^{N_0}_{N_0/p}\Wc_{p'}\Goo^\#\rp = 0,
    \end{align*}
    where $N_2' := N_2/p'$ and the last equality follows from $\Goo$, hence $\Wc_{p'}\Goo^\#$, being a newform of level $N_0$.
    Therefore, we can suppose $N_0 \mid N_1$, and obtain
    $$
    I _{\Delta, N_1, N_2, N_3}(g', \Goo^\#)
    =
      I _{\Delta, N_1, N_2, N_3}(g', \Wc_{N_0}\Goo^\#)    
    $$
    as a consequence of \eqref{eq:AL-op}.

    Set $N_3 = N_{1, 3} N_{2, 3} N_3'$ with $N_{j, 3} := \gcd(N_j, N_3)$. 
Using $\gcd(N_0, N_3) = 1$ and applying \eqref{eq:Vd} gives us
    \begin{align*}
      I_{\Delta, N_1, N_2, N_3}
      &(g_\tau', \Goo^\#)
=
        I_{\Delta, N_1, N_2, N_3}(g_\tau', \rho(w_{N_{1, 3}})\Goo^\#)\\
&=  N_3^{-k}      I_{\Delta, N_1, N_2, 1}
  \lp g_\tau', \rho\lp \smat{1/N_3}{}{}1_{N_3}  w_{N_{1,3}}\rp  \Goo^\#\rp\\
&=
N_3^{-k}  \frac{\vol(\Gamma_0(N)\backslash \Hb)}{\vol([H_0])}
\int_{[H_0]} \Goo^\#(h)
\overline{  \theta_{V_0}\lp g_\tau', h, \omega_{0}\lp w_{N_{1,3}} \smat{N_3}{}{}1_{N_3}  \rp  \varphi^k_{  \Delta, N_1, N_2}\rp}dh.
    \end{align*}
It is easy to check that
$$
\omega_{0}\lp w_{N_{1,3}} \smat{N_3}{}{}1_{N_3}  \rp
\varphi^\infty_{  \Delta, N_1, N_2}
=
\varphi^\infty_{  \Delta, N_1, N_2'}
$$
with $N_2' := N_2/N_{2, 3}$. 
Then $N' = \mathrm{lcm}(N_0, N_1N_2') = N/N_{2,3}$ and we have
\begin{equation}
  \label{eq:I-DN12}
      I_{\Delta, N_1, N_2, N_3}      (g_\tau', \Goo^\#)
= N_3^{-k} \sigma_1(N_{2, 3})            I_{\Delta, N_1, N_2', 1}      (g_\tau', \Goo^\#). 
\end{equation}
This takes care of $N_3$ and we omit it from the notations.

    To handle $N_2$, we use \eqref{eq:tr-Hecke} and proceed as follows
    \begin{align*}
      I _{\Delta, N_1, N_2}(g', \Goo^\#)
      &=
\frac{\vol(\Gamma_0(N_1)\backslash \Hb)}{\vol([H_0])}
\int_{[H_0]} \Goo^\#(h)
        \overline{\tr^{N}_{N_1}  \theta_{V_0}(g_\tau', h, \varphi^k_{  \Delta, N_1, N_2})}dh\\
      &=
\frac{\vol(\Gamma_0(N_1)\backslash \Hb)}{\vol([H_0])}
\int_{[H_0]} \Goo^\#(h)
        \overline{\Tc_{N_2}  \theta_{V_0}(g_\tau', h, \varphi^k_{  \Delta, N_1, 1})}dh\\
      &=       I _{\Delta, N_1, 1}(g', \Tc_{N_2} \Goo^\#)
        =   N_2^{1-k}a_G(N_2)   I _{\Delta, N_1, 1}(g', \Goo^\#). 
    \end{align*}
So we can suppose $N_2 = 1, N_0 \mid N= N_1$ and omit $N_2$ from the notations. 

For any prime $p \mid N/N_0$, we can compute on the Schwartz function $\varphi^\infty_{  \Delta, N}$  at $p$ to show that
\begin{equation}
  \label{eq:trace-id}
  \tr^{N}_{N/p} \varphi_{  \Delta, N}^\infty
  =
  \Tc_p \varphi_{  \Delta, N/p}^\infty
  - p^{3/2} \Uc_{p^2} \varphi_{  \Delta, N/p}^\infty
  + \varphi_{  \Delta, N/p}^\infty. 
\end{equation}
This gives us
\begin{align*}
        I _{\Delta, N}(g', \Goo^\#)
      &=
\frac{\vol(\Gamma_0(N_0)\backslash \Hb)}{\vol([H_0])}
\int_{[H_0]} \Goo^\#(h)
        \overline{\tr^{N}_{N_0}  \theta_{V_0}(g_\tau', h, \varphi^k_{  \Delta, N})}dh\\
      &=
        \sum_{c \mid (N/N_0)}\sum_{\rr \mid c}
(-1)^{\omega(\rr)}
\rr^{3/2}        (        \Uc_{\rr^2} I _{\Delta, N_0})(g',  \Tc_{N/(N_0 c)}\Goo^\#)        ,
\end{align*}
and leaves us with the case $N = N_1 = N_0, N_2 = 1$. 
For this, we apply Theorem 6.1 in \cite{AS21} to obtain the $(m/4)$-th Fourier coefficient of $\sqrt{v}^{-k+1/2}I_{\Delta, N}$ as
  \begin{align*}
&      (-1)^{k-1} \sqrt{|\Delta|}^{1-k}
   \sum_{\mu \in (\Lo^\vee \cap \frac1{2|\Delta|}\Lo)/\Lo}
   \chi_\Delta(\mu)
   \sum_{X \in \Gamma_0(N) \backslash \Lo + \mu,~ 4Q_0(X) = -m}
  \int_{c(X)}\Goo(z) \lp X, \smat{z^2}zz1 \rp^{k-1} dz\\
&=     (-1)^{k-1} \sqrt{|\Delta|}^{1-k}
  \sum_{\substack{\Qf = [Nb, -c, a]\in \Gamma_0(N) \backslash \Zb^3\\
    \mathrm{disc}(\Qf) = |\Delta|m,~ N \mid \Qf(1, 0)     }}
    \chi_\Delta(\mu)
    \int_{c(\Qf)}\Goo(z)
    (Nb z^2 - cz + a)^{k-1} dz\\
&=
  (-1)^{k-1+\lfloor k/2 \rfloor}
  \sqrt{|\Delta|}^{1-k}
        2^{-k} 
  \frac{|\Goo|_{\mathrm{Pet}}^2}{6|\tGoo|_{\mathrm{Pet}}^2}
\overline{  c_{\tGoo}(|\Delta|)} c_{\tGoo}(m),
  \end{align*}
  Here $X = \frac1{|\Delta|} \smat{a}{c/2}{c/2}{Nb}$ corresponds to the binary quadratic form $\Qf = [Nb, -c, a]$, and
  $-m = 4Q_0(X) = -|\Delta|(c^2 - 4Nab)$ implies that $m \in \Nb$ and $(-1)^k m \equiv 0, 1 \bmod 4$. 
  Also, the last step is a result of Kohnen \cite[Theorem 3]{Kohnen-FC}.
  The factor of 6 comes from different normalization of $|\tGoo|_{\mathrm{Pet}}^2$. 
Putting these together and applying equations \eqref{eq:Tcp} and \eqref{eq:Ucm} finishes the proof.
\end{proof}

\begin{remark}
  \label{rmk:ep}
  Suppose $2 \nmid \Delta$. For $\ep = 0, 1$, we define
  \begin{equation}
    \label{eq:varphi-ep}
    \varphi_{  \Delta, N_1, N_2}^{\infty, \ep}
    :=
    \varphi_{  \Delta, N_1, N_2}^{\infty}
    \cdot
            \cha\lp\Syme(\tfrac1{2\Delta N_1N_2} \hat\Zb) \cap \smat{N_2^{-1}\hat\Zb}{\hat\Zb + \frac{\ep}2}{\hat\Zb + \frac{\ep}2}{N_1N_2\hat\Zb} \rp.
  \end{equation}
  Then $    \varphi_{  \Delta, N_1, N_2}^{\infty} =     \varphi_{  \Delta, N_1, N_2}^{\infty, 0} +     \varphi_{  \Delta, N_1, N_2}^{\infty, 1}$. Using this Schwartz function, we can define $\theta^{k, \ep}_{\Delta, N_1, N_2}$ and $I^\ep_{\Delta, N_1, N_2, N_3}$ as in \eqref{eq:theta-Delta} and \eqref{eq:IDel}.
  If $\theta^{(k)}_{\Delta, N_1, N_2}$ on the left hand side of \eqref{eq:Shintani-eigenform} is replaced by $\theta^{k, \ep}_{\Delta, N_1, N_2}$, then $U_{\rr^2}\tGoo$ on the right hand side will be replaced by $U_{\rr^2}\tGoo^\ep$, where
  \begin{equation}
    \label{eq:Goo-ep}
    \tGoo^\ep(\tau) :=
\frac12 \lp     \tGoo(\tau)
+ (-1)^\ep    \tGoo(\tau + \tfrac12) \rp
= \sum_{n \in \Zb,~ n \equiv \ep \bmod 2} c_{\tGoo}(n)q^n. 
  \end{equation}

\end{remark}

\subsection{Generalized Rankin-Selberg L-series}


  We start by describing the action of the Cohen operators on the kernel $\Cc_{k-1}(\theta(g, h, \varphi^{(0, 1)})$.
In addition to $V_0$ from \eqref{eq:V0}, we need the quadratic spaces
\begin{equation}
  \label{eq:V1}
V_{1}=(\mathbb{Q},|D_2|\cdot x^2)
\end{equation}
with Weil representation $\omega_1 = \omega_{V_1}$. 
For $\ell \in \mathbb{N}_0$ and $\varphi_1 \in \mathcal{S}(\hat{V}_{1})$, we have the definite theta function
\begin{equation}
  \label{eq:theta1l}
  \begin{split}
    \theta_{1, D_2}^{(\ell)}\left(\tau, \varphi_1\right)
& :=
(2 \sqrt{v})^{-\ell} \sum_{m \in V_{1}(\mathbb{Q})} \varphi_1(m) \mathrm{H}_{\ell}\left(\sqrt{4 |D_2| \pi v}m\right) q^{|D_2|m^2}\\
&    = (-R_\tau)^{\lfloor \ell /2 \rfloor}
  \theta_{1, D_2}^{(\ell \bmod 2)}\left(\tau, \varphi_1\right)
  =
  \sqrt{v}^{-(\ell + 1/2)}\theta_{V_1}(g_\tau', \varphi_1 \otimes   \varphi^\ell_{1,\infty})  ,    
  \end{split}
\end{equation}
where
$\varphi_{1, \infty}^{\ell}(x)=2^{-\ell} \mathrm{H}_{\ell}(2\sqrt{\pi} x) e^{-2\pi x^2} \in
\mathcal{S}(V_{1}(\mathbb{R}))$, $g_\tau' = [g_\tau, 1] \in G'_\Rb$
and
$$\mathrm{H}_{\ell}(\xi):=(-1)^{\ell} e^{\xi^2}\left(\frac{d}{d \xi}\right)^{\ell} e^{-\xi^2}=\left(2 \xi-\frac{d}{d \xi}\right)^{\ell} \cdot 1$$
is the $\ell$-th (physicist's) Hermite polynomial.
In addition, we define
\begin{equation}
  \label{eq:theta1ep}
  \theta^{\ell, \ep}_{1, D_2}(\tau) := \theta_{1, D_2}^{(\ell)}(\tau, \cha( \hat \Zb + \tfrac{\ep}2)),~
  \theta^{\ell, \ep}_{V_1}(g_\tau') :=
 \sqrt{v}^{\ell + 1/2}    \theta^{\ell, \ep}_{1, D_2}(\tau) 
\end{equation}
for $\ep = 0, 1$. 

Set $d=|\Delta|=|D_1|$. Using the decomposition
\begin{equation}
  \label{eq:V-decomp}
  V \cong V_0 \oplus V_1,~
  (a, b, \nu)
  \mapsto
  \lp
  \left(\begin{array}{cc} a&\frac{\tr(\nu)}{2} \\ \frac{\tr(\nu)}{2} &b\end{array}\right)
, \frac{\nu-\nu^\prime}{2\sqrt{D_2/D_1}}
\rp,
\end{equation}
we identify $\Sc(V(\Ab)) \cong \Sc(V_0(\Ab)) \otimes \Sc(V_1(\Ab))$. 
By computation in the Fock model of $\omega_V$, we have the following result.

\begin{lemma}
  \label{lemma:RC-decomp}
Given $\varphi_j \in \Sc(\hat V_j)$ with $j = 0, 1$ and odd $k \in \Nb$, we have
  \begin{equation}
    \label{eq:RC-decomp}
    (  \Cc_{k-1} \theta)(g, h_0, (\varphi_0 \otimes \varphi_1) \otimes\varphi_\infty^{0, 1})
    -
\alpha_k
    \theta^{(k)}_{V_0}(g, h_0, \varphi_0)
    \theta_{V_1}(g, \varphi_1 \otimes \varphi_{1, \infty}^{k-1}    )  \in \mathrm{Im}(R^{H_0}),
  \end{equation}
    where
$\alpha_k := (-4\pi)^{(1-k)/2}\binom{ - k/2}{(k-1)/2 } $
and    $R^{H_0}$ is the raising operator on $H_0(\Rb) \cong \PGL_2(\Rb)$. 
\end{lemma}

\begin{proof}
  We can prove this using the Fock model in section 4.1 \cite{BLY22}. In the notation loc.\ cit.,
  $\iota(\Cc_{k-1} \varphi_\infty^{0, 1})$ is
  $$
  -\frac{\sqrt2 i}{4\pi}
 (-1)^{(k-1)/2}
  (32\pi^2)^{1-k}   \wf^k
 \sum_{s = 0}^{(k-1)/2} \binom{ - k/2}{(k-1)/2 - s} \binom{k/2 - 1}{s} (\vf - \bar\vf)^{k-1-2s}(\vf + \bar\vf)^{2s}  .
 $$
 The terms with $s > 0$ are in $ \mathrm{Im}(R^{H_0})$.
 Also, $\iota(\varphi_{0, \infty}^k) = (2\sqrt2 \pi i )^{-k} \wf^k$ (see e.g.\ (2.50) in \cite{BLS24}
 \footnote{Note our $\varphi^{k}_{0, \infty}$ is $(-\sqrt2)^k$ times $\varphi^{(0, k)}_\infty$ in \cite{BLS24}.}), and $(8\pi)^{(1-k)/2} (i\zf_2)^{k-1} = (-R)^{(k-1)/2} 1$ in the Fock model.
 Putting these together finishes the proof.
\end{proof}

With \eqref{eq:tMv}, this gives us
\begin{equation}
  \label{eq:tV-decomp}
  \tV \cong \k_2 \oplus V_0 \oplus V_1,~
  (\alpha, a, b, c, d) \mapsto
  \lp
  \frac{\alpha}{ND_2},
\frac1{D_1} \pmat{N^{-1}a}{c/2}{c/2}{b}, \frac{d}{2D_2}\rp.
\end{equation}
We can write $\tv^\infty_{N } = \sum_j\phi_{2, N, j}\otimes  \phi_{N, j}$ with $\phi_{2, N, j} \in \Sc(\hat \k_2), \phi_{N, j} \in \Sc(\hat V)$.
This decomposition is not canonical and could be complicated. Fortunately, equation \eqref{eq:tvp-val} tells us that 
\begin{equation}
  \label{eq:tvp-decomp-val}
  \begin{split}
      \tv^\infty_{N} \mid_{(\nf_2 \df_2^{-1} \oplus V)\otimes \hat\Zb}
    &= \phi_{N} :=
      \sum_{\substack{\ep_0, \ep_1 \in \{0, 1\}\\ \frac2{\gcd(2, D_2)} \mid(\ep_0 + \ep_1) }}
        \varphi^{\infty, \ep_0}_{  D_1, N, 1} \otimes \cha(\hat\Zb + \tfrac{\ep_1}2),
  \end{split}
\end{equation}
where $\varphi^{\infty, \ep}_{  \Delta, N_1, N_2}$ is defined in \eqref{eq:varphi-ep} and Remark \ref{rmk:ep}.

Now, we are ready to evaluate the $L$-function $L(s, \Go; D_1, D_2)$
defined in \eqref{eq:LsG}.
First, using equation \eqref{eq:matchN} with $k = 1$ gives us
\begin{equation}
  \label{eq:LsG1}
    L(s, \Go; D_1, D_2)
 =
          \sum_{N_1 \mid N}
\frac{4 {\zeta_N(s+1)}
  {\beta_{N/N_1}(s)}}{ {\zeta_N(1)}
    \Lambda^{}(s+1, \chi_1)
}
          \int_{\Gamma_0(N)\backslash \Hb} \Go(z)
   \overline{\Cc_{k-1}(
I(h_z, \tv^1_{N_1}, \bar s))} y^{k}d\mu(z) ,
\end{equation}
where $\beta_{ N'}(s) :=
\frac{ \zeta_{N'}(1)}{ N'\zeta_{N'}(s)}$ and $\Lambda(s, \chi_1)$ is the completed $L$-function defined below equation \eqref{eq:E*}. 
For each $N_1 \mid N$, we apply the decomposition of $\tv_{N_1}^\infty$ below \eqref{eq:tV-decomp}, its restriction in \eqref{eq:tvp-decomp-val}, interchange order of integration, unfold and apply Lemma \ref{lemma:RC-decomp} to obtain
\begin{align*}
&  \int_{\Gamma_0(N)\backslash \Hb}
  \Go(z)
  \overline{\Cc_{k-1}(
  I(h_z, \tv^\infty_{N_1}, \bar s))} y^{k}d\mu(z)
  \\
&=  \int_{\Gamma_0(N)\backslash \Hb} \Go(z)
   \overline{
  \int_{[\SL_2]}
\sum_j
  E(g, \lambda(\phi_{2, N_1, j}) \otimes \Phi_\infty^1, \bar s)
  \Cc_{k-1}(  \theta)(g, h_z, \phi_{N_1, j} \otimes \varphi_\infty^{(0, 1)})dg } y^{k}d\mu(z)\\
&=
  \int_{[B]/B(\hat\Zb)}
  \sum_j
  \overline{
  \omega(g)(\phi_{2, N_1, j})(0) \otimes \Phi_\infty^1(g  , \bar s)}
  \int_{\Gamma_0(N)\backslash \Hb} \Go(z)
  \overline{
  \Cc_{k-1}(  \theta)(g, h_z, \phi_{N_1, j} \otimes \varphi_\infty^{(0, 1)}) } y^{k}d\mu(z) dg \\ 
&=
  \int_{B(\Zb)\backslash B(\Rb)}
  {|a(g)|^{1+ s}}
  \int_{\Gamma_0(N)\backslash \Hb} \Go(z)
  \overline{
 \alpha_k \theta_{V}(g, h_z, \phi_{N_1} \otimes( \varphi_{0, \infty}^{k}\varphi_{1, \infty}^{k-1}))} y^{k}d\mu(z) dg.
\end{align*}
We can set $g = g_\tau$ for $\tau = u + iv \in \Hb$. Then $a(g) = \sqrt{v}, dg = \frac{d\mu(\tau)}2$ and $B(\Zb)\backslash B(\Rb) = [0, 1] \times \Rb_{> 0 }$. 

Suppose $\Go(z) = \Goo(N_3z)$ for a newform $\Goo \in S_{2k}(N_0)$ and some $N_3 \mid (N/N_0)$. 
Setting $g = g_\tau$, applying \eqref{eq:theta-Delta1} and the definition of $\phi_{N_1}$ in \eqref{eq:tvp-decomp-val}, we can rewrite the integral over $\Gamma_0(N)\backslash \Hb$ as
\begin{align*}
  \int_{\Gamma_0(N)\backslash \Hb}
  &\Go(z)
  \overline{
 \theta_{V}(g_\tau, h_z, \phi_{N_1} \otimes( \varphi_{0, \infty}^{k}\varphi_{1, \infty}^{k-1}))} y^{k}d\mu(z) \\
&=
\sigma_1(N/N'')
  \sum_{\substack{\ep_0, \ep_1 \in \{0, 1\}\\ \frac2{\gcd(2, D_2)} \mid(\ep_0 + \ep_1) }}
  I^{\ep_0}_{D_1, N_1, 1, N_3}(g'_\tau, \Goo)
  \overline{\theta^{k-1,\ep_1}_{V_1}(g'_\tau)}\\
&=
  \sqrt{v}^{-1}
  \sigma_1(N/N'')
  \sum_{\substack{\ep_0, \ep_1 \in \{0, 1\}\\ \frac2{\gcd(2, D_2)} \mid(\ep_0 + \ep_1) }}
      \int_{\Gamma_0(N'') \backslash \Hb} \Goo(N_3 z)
      \overline{\theta^{k, \ep_0}_{D_1, N_1, 1}(\tau, z     )} y^{2k}d\mu(z)
\overline{  \theta^{k-1,\ep_1}_{1, D_2}(\tau)},
\end{align*}
where $N'' := \mathrm{lcm}(N_0N_3, N_1)\mid N$, 
$I^\ep_{\Delta, N_1, N_2, N_3}$ is defined in \eqref{eq:IDel} and Remark \ref{rmk:ep}, and $\theta_{1, D_2}^{k-1, \ep_1}$ is defined in \eqref{eq:theta1l} and \eqref{eq:theta1ep}. 
The factor $\sigma_1(N/N'')$ is the index of $\Gamma_0(N)$ in $\Gamma_0(N'')$.
By Proposition \ref{prop:Shintani}, this integral does not vanish only if $N_0 \mid N_1$, in which case
\begin{equation}
  \label{eq:LsG2}
  \begin{split}
      &  \int_{\Gamma_0(N'')\backslash \Hb}
  \Go(z)
  \overline{\Cc_{k-1}(
  I(h_z, \tv^1_{N_1}, \bar s))} y^{k}d\mu(z)
    \\
      &=
\alpha_k        D_1^{(1-k)/2}
\overline{c_{\tGoo}(|D_1|)}
        2^{-k} 
        \frac{|\Goo|_{\mathrm{Pet}}^2}{6|\tGoo|_{\mathrm{Pet}}^2}
N_3^{-k}
        \lp\frac{N''}{N_0}\rp^{1-k}
        \sum_{c \mid (N/N_0)}
        (-1)^{\omega(\rr)}        a_{\Goo}\lp \frac{N''}{N_0 c}\rp\\
&\times        \int_0^\infty
\sqrt{v}^{2k-1 + s}
   \sum_{\substack{\ep_0, \ep_1 \in \{0, 1\}\\ \frac2{\gcd(2, D_2)} \mid(\ep_0 + \ep_1) }}
\int_0^1
        \sum_{\rr \mid c}
    \lp U_{\rr^2}        \tGoo^{\ep_0}\rp(\tau/4)
\overline{  \theta^{k-1,\ep_1}_{1, D_2}(\tau)}
        \frac{du dv}{2v^2}
  \end{split}
\end{equation}
Note $N' = N$ as $N_2 = 1$ and we have used $2 \nmid k$. 
The last line above simplifies as
\begin{align*}
&  \int_0^\infty
\sqrt{v}^{2k-1 + s}
   \sum_{\substack{\ep_0, \ep_1 \in \{0, 1\}\\ \frac2{\gcd(2, D_2)} \mid(\ep_0 + \ep_1) }}
\int_0^1
    \lp U_{\rr^2}        \tGoo^{\ep_0}\rp(\tau/4)
\overline{  \theta^{k-1,\ep_1}_{1, D_2}(\tau)}
  \frac{du dv}{2v^2}\\
  &= 
(-4)^{(1-k)/2}  \int_0^\infty
\sqrt{v}^{k +    s}
2 \sum_{n \ge 1} c_{\tGoo}(d^2n^2|D_2|) H_{k-1}(\sqrt{|D_2|\pi v} n) e^{-\pi v |D_2| n^2}
  \frac{ dv}{2v}\\
  &= 
    (-4)^{(1-k)/2} \sqrt{|D_2|\pi }^{-k-   s}
2 \sum_{n \ge 1} \frac{c_{\tGoo}(d^2n^2|D_2|)}{n^{k+   s}}
    \int_0^\infty
    t^{k +    s}
    H_{k-1}(t) e^{-t^2}
  \frac{ dt}{t}.
\end{align*}
Using Theorem 2 in \cite{Kohnen-newform} and $d \mid N$ satisfies \eqref{eq:Heegner}, we can deduce that
\begin{equation}
  \label{eq:Kohnen-new}
  \sum_{n \ge 1} \frac{c_{\tGoo}(\rr^2n^2|D_2|)}{n^{k+s}}
  =
  \frac{c_{\tGoo}(|D_2|)}{L(   s + 1, \chi_2)} L(   s + k, \Goo)
\delta_\rr(\Goo, s),
\end{equation}
where $\delta_\rr(\Goo, s) := \zeta_\rr(   s + 1) \prod_{p \mid \rr \text{ prime}} (a_{\Goo}(p) - p^{k-1} (1+p^{-s}))$ is multiplicative. 
Also, it is an easy exercise to show that $    \int_0^\infty
    t^{k + s}
    H_{k-1}(t) e^{-t^2}
    \frac{ dt}{t} = \frac{\Gamma(k+s)}{\Gamma(1+s)}$.
    Putting these together gives us the following result.

    \begin{theorem}
      \label{thm:LsG}
      Let $D_1, D_2 < 0$ be distinct fundamental discriminants, not both even.
      For $G \in S_{2k}(N)$ with  $N$ square-free and satisfying \eqref{eq:Heegner}, 
let      $L(s, \Go; D_1, D_2)$ be the $L$-function defined in \eqref{eq:LsG}. 
      Suppose $\Go(z) = \Goo(N_3z)$ for a newform $\Goo \in S_{2k}(N_0)$ and some $N_0 \mid N$ and $N_3 \mid (N/N_0)$.
      Then
      \begin{equation}
        \label{eq:LsG0}
        \begin{split}
                  L(s, \Go; D_1, D_2)
 &       =
\frac{|\Goo|_{\mathrm{Pet}}^2}{3|\tGoo|_{\mathrm{Pet}}^2}
   \frac{(1 + \epsilon(\Goo))\overline{c_{\tGoo}(|D_1|)}{c_{\tGoo}(|D_2|)}}{\Lambda(s+1, \chi_1)\Lambda(s+1, \chi_2 )|D_1D_2|^{(k-1)/2}
   }
   L(\Goo, s + k)
          \\
 &\quad
\times   C_k(s) 
   \frac{ {\zeta_N(s+1)}}{ {\zeta_N(1)}}
   \prod_{p \mid (N/(N_0 N_3))} \gamma_{p, 0}(\Goo, s)
   \prod_{p \mid N_3} \gamma_{p, 1}(\Goo, s)      
,
        \end{split}
\end{equation}
where
$  C_k(s) :=
 \binom{-k/2}{(k-1)/2} 2^{3-3k} \pi^{-k-s-1/2} \frac{\Gamma(k+s)}{\Gamma(1+s)}
\Gamma(s/2 + 1)$ and 
\begin{equation}
  \label{eq:gammas}
  \begin{split}
\gamma_{p, 0}(\Goo, s)
&   := 
    \frac{-p^{-k-s}}{1-p^{-1-s}} a_{\Goo}(p) + \frac{(1 + p^{-1})p^{-2s} - (3 + p^{-1})p^{-s} + 2p}{(p-1)(1-p^{-1-s})},\\
\gamma_{p, 0}(\Goo, s)
&   := 
    p^{-2k}\lp \frac{p^{-2s} - 2p^{1-s} + p}{(p-1)(1-p^{-1-s})} a_{\Goo}(p) + p^k \frac{1 + p^{-s}}{1 - p^{-1-s}} \rp. 
  \end{split}
\end{equation}
    \end{theorem}
    \begin{proof}
      Following the calculations before the theorem, we arrive at
\begin{align*}
L(s, \Go; D_1, D_2) =& \frac{|\Goo|_{\mathrm{Pet}}^2}{3|\tGoo|_{\mathrm{Pet}}^2}
\frac{(1 + \epsilon(\Goo))\overline{c_{\tGoo}(|D_1|)}{c_{\tGoo}(|D_2|)}}{D^{\frac{k-1}{2}} \Lambda(s+1, \chi_1)\Lambda(s+1, \chi_2)}
L(\Goo, s + k)\\
&\times C_k(s) \frac{ {\zeta_N(s+1)} }{ {\zeta_N(1)}} B_{N/N_0}(\Goo, N_3, s), 
\end{align*}
where we have
$\beta_{N'}$ defined as in \eqref{eq:LsG1} and
\begin{equation}
  \label{eq:BN}
  \begin{split}
      B_{N'}(\Goo, N_3, s) :=
             N_3^{-k}
&  \sum_{N_1' \mid N'} \sigma_1(N'/\mathrm{lcm}(N_1', N_3))
  {\beta_{N'/N^\prime_1}(s)}
        \mathrm{lcm}(N_1', N_3)^{1-k}\\
            &\times  \sum_{c \mid {N_1'}{}}
              a_{\Goo}\lp \frac{\mathrm{lcm}(N'_1, N_3)}{ c}\rp
              \sum_{\rr \mid c}
              (-1)^{\omega(\rr)}
              \delta_\rr(\Goo, s).
    \end{split}
\end{equation}
Since $N$ is square-free and $N_3 \mid N' = N/N_0$,
we can  write $N_1', c, \rr$ as products of factors co-prime to  $N_3$ and dividing $N_3$.
This gives us $B_{N'}(\Goo, N_3, s) = B_{N'/N_3, 0}(\Goo, s) B_{N_3, 1}(\Goo, s)$, where
\begin{align*}
  B_{M, 0}(\Goo, s)
  &:=
    \sum_{M_1 \mid M}
    \sigma_1(M/M_1)
    {\beta_{M/M_1}(s)}
    M_1^{1-k}
    \sum_{c \mid M_1}
    a_{\Goo}\lp \frac{M_1}{ c}\rp
    \sum_{\rr \mid c}
    (-1)^{\omega(\rr)}
\delta_\rr(\Goo, s),\\
  B_{M, 1}(\Goo, s)
  &:=
        M^{-k}
    \sum_{M_1 \mid M}
  {\beta_{M/M_1}(s)}
        M^{1-k}
    \sum_{c \mid {M_1}{}}
  a_{\Goo}\lp \frac{M}{ c}\rp
    \sum_{\rr \mid c}
    (-1)^{\omega(\rr)}
    \delta_\rr(\Goo, s).
\end{align*}
Both functions $B_{M, 0}, B_{M, 1}$ are multiplicative in $M$, and it is straightforward to check that $\gamma_{p, j}(\Goo, s) = B_{p, j}(\Goo, s)$. This completes the proof. 
    \end{proof}

\subsection{Proof of Theorem \ref{thm:sing}}
Suppose $\Psi = \Psi_f$ with $f=\sum_{m \gg -\infty}c_f(m) q^{m} \in M^!_0(N)$ such that it has pole only at $\infty$, satisfies $c_f(m) \ge 0$ for all $m < 0$, and $\vv(f)$ has trivial constant term at the trivial coset. 
By Theorem \ref{thm:CMformula} with $k=1$ and Equation \eqref{eq:CM-val}, we have
\begin{align*}
\log | \Nm_{H/\Qb} \Psi_f(\tau_1, \tau_2)|
=
\log |\Psi_f(Z(W))|
=
  -\frac{|Z(W)|}{2}
  \mathrm{CT}\left(\tr^N_1\lp f \cdot\Ec_N^{+}\rp\rp. 
\end{align*}
Since $\vv(f)$ has no constant term at the trivial component and the vector-valued incoherent Eisenstein series $E'_\Nf(z, 0)$ has no constant term at all non-trivial components \cite[Proposition 4.6]{BKY12}, we can apply  \eqref{eq:CTvv} to see that the summation on the right hand side involve only Fourier coefficients of $\Ec^+_N$ (at various cusps) with positive index. 


To show that this is non-zero, it suffices to show that
$  \mathrm{CT}\left(\tr^N_1\lp \tilde f \cdot\Ec_N^{+}\rp\rp \neq 0$
  by the same argument in section 4 of \cite{Li21}, where $\tilde f \in H_{-4}(N)$ is the unique harmonic Maass form having the same principal part as $f$. 
  For $N = 1$ in \cite{Li21}, the form $\tilde f$ is weakly holomorphic. For $N > 1$, it is in general not weakly holomorphic.

  Again by Theorem \ref{thm:CMformula} with $k=3$, we have
  \begin{align*}
      -\frac{2}{|Z(W)|}
G_{k, \tilde f}(Z(W))
=
\mathrm{CT}\left(\tr^N_1\lp \tilde f^{+} \cdot\mathcal{C}_{k-1}\left(\Ec_N^{+}\rp\rp\rp
    - L'(0, \xi(\tilde f); D_1, D_2). 
  \end{align*}
  Using the positivity of $G_{k, \tilde f}$ and equidistribution of $Z(W)$ on $X_0(N)^2$ as $\max(|D_1|, |D_2|)$ goes to infinity, we can find absolute constant $C > 0$ such that the absolute value of the left hand side is bounded below by $C$.
  On the other hand, since $\xi(\tilde f) \in S_6(N)$ is independent of $D_i$, we can apply Theorem \ref{thm:LsG}, Brauer-Siegel Theorem and non-trivial bounds for Fourier coefficients of half-integral weight eigenforms \cite{Duke1988} to obtain
  \begin{align*}
    | L'(0, \xi(\tilde f); D_1, D_2)|
    &\ll_N \sum_{\Go \text{ eigenform in } S_6(N)} |L'(0, \Go; D_1, D_2)|\\
    &\ll_N
      \sum_{N_0 \mid N,~ \Goo \text{ newform in } S_6(N_0)}
\frac{1}{|\tGoo|_{\mathrm{Pet}}^2}
      \left|
      \frac{c_{\tGoo}(|D_1|)c_{\tGoo}(|D_2|)}{\Lambda(1, \chi_1)\Lambda(1, \chi_2) |D_1D_2|^{(k-1)/2}}
      \right|\\
    &\ll_N |D_1D_2|^{-\delta}
  \end{align*}
  for some fixed $\delta > 0$.
  This completes the proof.


\bibliography{DLIP.bib}{}
\bibliographystyle{amsplain}

\end{document}